\documentclass[12pt]{article}

\usepackage[english]{babel}
\usepackage{amsmath,amssymb,amsthm}
\usepackage{bm}
\usepackage{epsfig}

\newtheorem{theorem}{Theorem}[section]
\newtheorem{proposition}[theorem]{Proposition}
\newtheorem{lemma}[theorem]{Lemma}

\theoremstyle{definition}

\newtheorem*{remark}{Remark}

\newcommand{\Prob}{\mathbf{P}}
\newcommand{\E}{\mathbf{E}}
\newcommand{\n}{\mathbf{n}}
\newcommand{\w}{\mathbf{w}}
\newcommand{\x}{\mathbf{x}}
\newcommand{\y}{\mathbf{y}}
\newcommand{\z}{\mathbf{z}}
\newcommand{\Z}{\mathbb{Z}} 
\newcommand{\R}{\mathbb{R}}
\newcommand{\C}{\mathbb{C}}
\newcommand{\Half}{\mathbb{H}}
\newcommand{\Disk}{\mathbb{D}}
\newcommand{\uplat}{\mathcal{H}}
\newcommand{\p}{\partial}
\newcommand{\dist}{\operatorname{dist}}
\newcommand{\excur}{\nu}
\newcommand{\eps}{\epsilon}

\newcommand{\hatgamma}{\hat{\boldsymbol{\gamma}}}
\newcommand{\SLE}{\text{SLE}}
\newcommand {\msae} {W}

\begin{document}

\title{The configurational measure on
mutually avoiding $\SLE$ paths}

\author{
Michael J.~Kozdron\footnote{Research supported in part by the Natural Sciences and Engineering Research Council of Canada.}\\
{\small Department of Mathematics \& Statistics, College West 307.14}\\
{\small University of Regina, Regina, SK S4S 0A2 Canada}\\ 
{\small \texttt{kozdron@math.uregina.ca}}\\        
\and
Gregory F. Lawler\footnote{Research supported by the National Science Foundation
grant DMS-0405021.}\\
{\small Department of Mathematics, 310 Mallot Hall}\\
{\small Cornell University, Ithaca, NY 14853 USA}\\
{\small \texttt{lawler@math.cornell.edu}}\\
}

\date{} 

\maketitle

\begin{abstract}
We define multiple chordal $\SLE$s in a simply connected domain by considering a natural configurational measure on paths. We show how to construct these measures so that they are conformally covariant and satisfy  certain boundary perturbation and Markov properties, as well as a cascade relation. As an example of our construction, we derive the scaling limit of Fomin's identity in the case of two paths directly; that is, we prove that the probability that an $\SLE_2$ and a Brownian excursion do not intersect can be given in terms of the determinant of the excursion hitting matrix. Finally, we define the $\lambda$-SAW, a one-parameter family of measures on self-avoiding walks on $\Z^2$.
\end{abstract}



\newpage


\section{Introduction}

The Schramm-Loewner evolution ($\SLE$) is a powerful tool for
describing continuum limits that arise in critical two-dimensional
models in statistical physics. The starting point~\cite{schramm}
is to write down some properties that the continuum limit is
expected to satisfy and then to show that this leaves only a one-parameter 
family of possibilities.  This works beautifully in
the case of one path in a simply connected domain.  

For multiple paths or non-simply connected domains, the 
approach in~\cite{schramm} is not sufficient to determine what
the natural definition(s) should be.  Essentially the conformal
Markov property,  which determines the distribution of 
a single $\SLE$ path in simply connected domains,  is not sufficient
to describe the more complicated situations.  A number of authors
have tried to consider these processes using infinitesimal descriptions.
In this paper, we would like to try to suggest that this is not
the most natural approach, and that instead one should look at configurational
measures on paths.  What we do here is describe the configurational
measure associated to multiple chordal $\SLE$s in a simply connected
domain.  This construction is based on the work in~\cite{restriction}
and is also discussed in~\cite{dubedat}. For a statistical mechanics argument describing multiple $\SLE$s, consult~\cite{BBK}. Nonetheless, it still seems useful to describe the configurational measure in some detail.

We will start by writing down a number of properties that we
expect the measure to have.  We will not worry at the moment whether
some of these properties are redundant or
consistent.   Let $D$ be a Jordan domain,
and let $z_1,\ldots,z_n, w_n, \ldots, w_1$
be $2n$ distinct
points on $\p D$ in counterclockwise order. We write 
$\z = (z_1,\ldots,z_n)$, $\w = (w_1,\ldots,w_n)$, and say that $D$ is (locally) 
analytic at $\z$, $\w$ if conformal
maps sending $D$ to the unit disk $\Disk$ can be extended analytically
 to neighborhoods
of $z_1,\ldots,z_n, w_n, \ldots,
w_1$.   We fix $b \in \R$, a parameter that is called the
{\em boundary scaling exponent} or {\em boundary conformal weight}.
We want to define a measure $Q_{D,b,n}(\z,\w)$ on 
mutually avoiding $n$-tuples of
simple paths $(\gamma^1,\ldots,\gamma^n)$ in $D$.    More
precisely, $\gamma^j$ is an equivalence class of curves such that
there is a representative $\gamma^j:[0,1] \rightarrow \C$ which
is simple, $\gamma^j(0) = z_j$, $\gamma^j(1) = w_j$, and $\gamma^j(0,1)
\subset D$.  We will write just $\gamma^j$ for $\gamma^j(0,1)$.
When we refer to measures on curves in this paper, we
will implicitly 
mean measures on equivalence classes of curves under
equivalence by reparametrization; consult, for example,~\cite{Kozdron} or~\cite{Lbook}.

\begin{figure}[htb]
\begin{center}
\epsfig{file=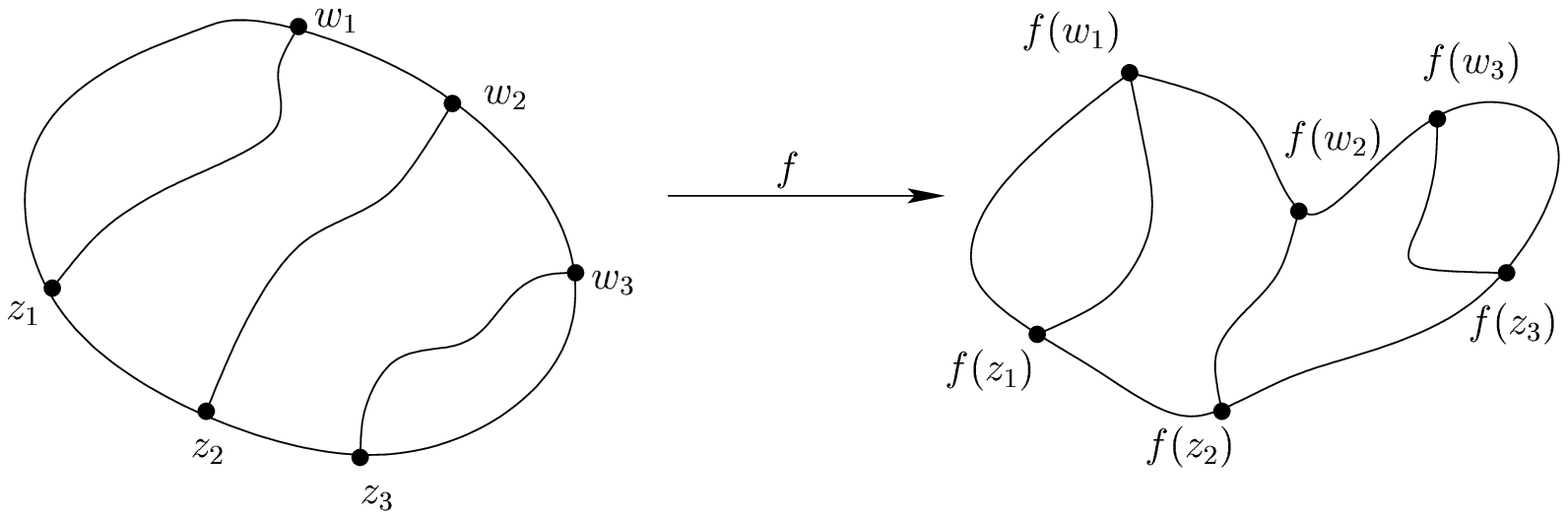}
\end{center}
\caption{Conformal Covariance}\label{figure1}
\end{figure}

\begin{itemize}
\item  \textbf{Conformal Covariance.}
If $D$ is analytic at $\z$, $\w$, then $Q_{D,b,n}(\z,\w)$
is a non-zero, finite measure supported on $n$-tuples
$(\gamma^1,\ldots,\gamma^n)$ where 
$\gamma^j$ is a simple curve in $D$ connecting $z_j$ and $w_j$ and
\[  \gamma^j \cap \gamma^k = \emptyset, \;\;\;\; 1 \leq j < k \leq n. \]
Moreover, if $f: D \rightarrow f(D)$ is a conformal transformation
and $f(D)$ is analytic at $f(\z)$, $f(\w)$, then
\begin{equation}  \label{confcov}
  f \circ Q_{D,b,n}(\z,\w) = |f'(\z)|^b \, |f'(\w)|^b
  \, Q_{f(D),b,n}(f(\z),f(\w))
\end{equation}
where $f(\z) = (f(z_1), \ldots, f(z_n))$ and $f'(\z) = f'(z_1) \cdots f'(z_n)$;
see Figure~\ref{figure1}.
\end{itemize}

If we write
\[     Q_{D,b,n}(\z,\w) = H_{D,b,n}(\z,\w) \, \mu^\#_{D,b,n}
   (\z,\w), \]
where $  H_{D,b,n}(\z,\w) = | Q_{D,b,n}(\z,\w)|$ and
$\mu^\#_{D,b,n} (\z,\w)$ is a probability measure, then the conformal
covariance condition~(\ref{confcov}) becomes the scaling rule for $H$,
\begin{equation}  \label{hinvariance}
        H_{D,b,n}(\z,\w) =   |f'(\z)|^b \, |f'(\w)|^b
  \,  H_{f(D),b,n}(f(\z),f(\w)), 
\end{equation}
and  the  conformal {\em invariance} rule for $\mu^\#$, namely
\begin{equation}  \label{invariance}
    f \circ \mu^\#_{D,b,n}(\z,\w)
   = \mu^\#_{f(D),b,n}(f(\z),f(\w))   . 
\end{equation}
Since $\mu^\#$ is a conformal invariant, we can  
  define  $\mu^\#_{D,b,n}(\z,\w)$ even if the boundaries
are not smooth at $\z$, $\w$. (Compare this with the construction of the Brownian excursion measure in~\cite{Kozdron}.)  
We will write $Q_{D,b}(z,w) = Q_{D,b,1}(z,w)$ and similarly
for $H_{D,b}(z,w)$ and $\mu^\#_{D,b}(z,w)$.

The next condition describes what happens when we remove
part of the domain.

\begin{figure}[htb]
\begin{center}
\epsfig{file=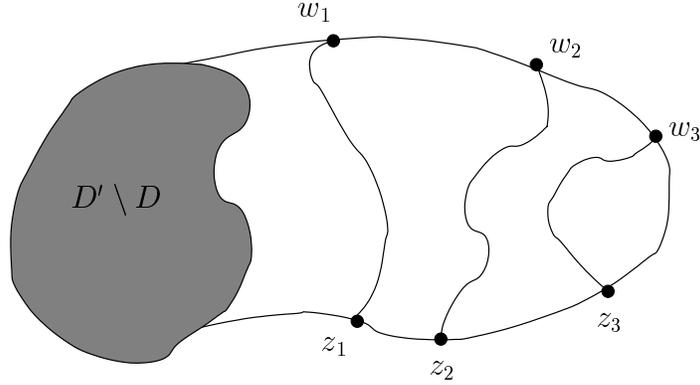}
\end{center}
\caption{Boundary Perturbation}\label{figure2}
\end{figure}

\begin{itemize}
\item  \textbf{Boundary Perturbation.} 
Suppose $D \subset D'$ are Jordan domains and
$\p D$, $\p D'$ agree and are analytic in neighborhoods
of $\z$, $\w$.  Then $Q_{D,b,n}(\z,\w)$ is absolutely
continuous with respect to $Q_{D',b,n}(\z,\w)$.  Moreover,
the Radon-Nikodym derivative
\[           Y_{D,D',b,n}(\z,\w)
 = \frac{dQ_{D,b,n}(\z,\w)}
              {dQ_{D',b,n}(\z,\w)} \]
is a conformal invariant. In other words, if $f: D' \rightarrow
f(D')$ is a conformal map that extends analytically in
neighborhoods of $\z,\w$, then
\begin{equation}\label{dec13eq1}
           Y_{f(D),f(D'),b,n}(f(\z),f(\w)) (f \circ
\bar \gamma) 
 = Y_{D,D',b,n}(\z,\w)(\bar \gamma), 
\end{equation}
where $\bar \gamma = (\gamma^1,\ldots,\gamma^n)$ and $f \circ \bar \gamma = (f \circ \gamma^1,\ldots,f \circ \gamma^n)$; see Figure~\ref{figure2}.
\end{itemize}

As with $\mu^\#_{D,b,n}(\z,\w)$,
the last condition~(\ref{dec13eq1}) implies that  $Y_{D,D',b,n}(\z,\w)$
is well-defined  even if the boundaries are not smooth at $\z,\w$.

The third relation relates the measures for different $n$. 

\begin{figure}[htb]
\begin{center}
\epsfig{file=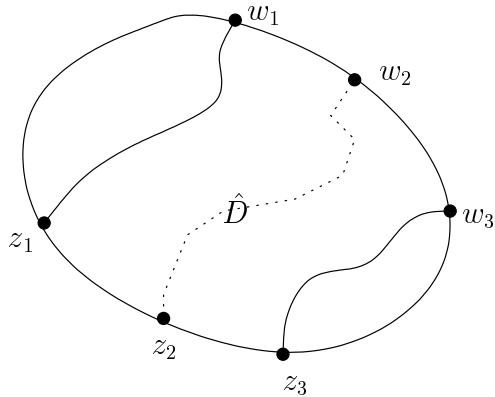}
\end{center}
\caption{Cascade Relation. $\hat z = (z_1,z_3)$, $\hat w = (w_1,w_3)$.}\label{figure3}
\end{figure}

\begin{itemize}
\item \textbf{Cascade Relation.} Let
 \[ \hat \z = (z_1,\ldots,z_{j-1}, z_{j+1},\ldots,
z_n) ,
\;\;\; \hat \w = (w_1,\ldots,w_{j-1},w_{j+1},\ldots,w_n),\]
\[   \hatgamma =  (\gamma^1,\ldots,\gamma^{j-1},
\gamma^{j+1},\ldots,\gamma^n). \]
The marginal distribution on  $\hatgamma$
 induced by  $Q_{D,b,n}(\z,\w)$
is absolutely continuous with respect to $Q_{D,b,n-1}(\hat \z,
\hat \w)$ with Radon-Nikodym derivative $H_{\hat D,b}(z_j,w_j)$.
Here $\hat D$ is the subdomain of $D \setminus \hatgamma$ 
whose boundary includes $z_j$, $w_j$; see Figure~\ref{figure3}. 
\end{itemize}

In order to construct such a collection of measures it is useful
to include a fourth condition.  This condition is   
stated in terms of  $\mu^\#_{D,b}$
rather than  $Q_{D,b}$.

\begin{figure}[htb]
\begin{center}
\epsfig{file=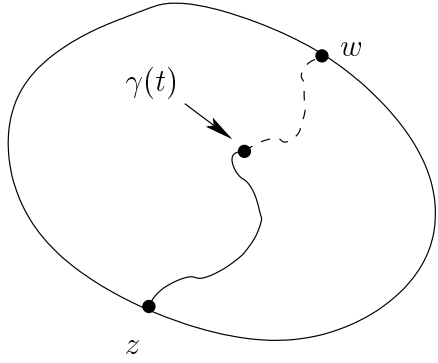}
\end{center}
\caption{Markov Property}\label{figure4}
\end{figure}

\begin{itemize}
\item \textbf{Markov Property.}  In the measure $\mu^\#_{D,b}(z,w)$,
the conditional distribution on $\gamma$ given an initial segment
$\gamma[0,t]$ is $\mu^\#_{D\setminus \gamma[0,t],b}(\gamma(t),w)$; see Figure~\ref{figure4}.
\end{itemize}

We have stated this condition in a way that does not use
two dimensions and conformal invariance.  The {\em conformal
Markov property} is the combination of the Markov property 
and~(\ref{invariance}). Schramm~\cite{schramm} showed that
there is a one-parameter family of measures, which he
parametrized by $\kappa$, satisfying the conformal Markov property.
While these measures are well-defined for $\kappa > 0$, they are supported
on simple curves only for $0 < \kappa \leq 4$; see~\cite[Proposition~6.9]{Lbook}. Therefore, we will restrict
our consideration only to $\kappa \leq 4$.  The scaling rule~(\ref{hinvariance}) 
implies that, up to an arbitrary multiplicative constant,
\begin{equation}  \label{3}
             H_{b}(x,y)=H_{\Half,b}(x,y) =H_{\Half,b,1}(x,y) =  |y-x|^{-2b}.
\end{equation}
We will continue to write $H_{b}$ for $H_{\Half,b}$ $(=H_{\Half,b,1})$ where $\Half = \{z \in \C:
\Im(z) > 0\}$ is  the upper half plane.  We will choose
the constant to be $1$ so that~(\ref{3}) holds.  Given this, the
measures $Q_{D,b,n}(\z,\w)$ can be determined for all simply
connected $D$ and $\z$, $\w$ at which $\p D$ is analytic by a
combination of conformal mapping and the cascade relation.  The
question is whether or not we can find $b$, $\kappa$ such that all
four of the properties above hold.  The answer is yes, assuming
we choose $\kappa = \kappa(b)$ satisfying
\[        \kappa = \frac{6}{2b+1} \;\;\; \text{so that} \;\;\;
    b = \frac{6 - \kappa}{2\kappa}. \]
Notice that $b$ decreases from $\infty$
to $1/4$ as $\kappa$ increases from $0$ to $4$ which 
puts the restriction $b \geq 1/4$ on the scaling exponent.
There are four other parameters that
 can  be used  to
parametrize this family, namely
\begin{equation}\label{paramdefns}
   a = \frac{2}{\kappa} = \frac{2b+1}{3}, \;\;\;
  \lambda = \frac{(3a-1) \, (4a-3)}{2a} , \;\;\; 
  c = - 2 \lambda, \;\;\; \text{and} \;\;\; d = 1 + \frac 1{4a} . 
\end{equation}
These parameters have interpretations which can be summarized as follows.
\begin{itemize}
\item   $b \in [1/4,\infty)$ is the \emph{boundary
scaling exponent} or \emph{boundary conformal weight}.
\item  $\kappa \in (0,4]$ is the variance of the driving Brownian motion in
   $\SLE_\kappa$ if the half-plane capacity at time $t$ is $2t$.
\item   $a \in [1/2, \infty)$ and $at$ is the half-plane capacity at time $t$ for $\SLE_\kappa$ if the driving Brownian motion is chosen to have variance $1$.
\item  $\lambda \in [-1/2,\infty)$ is the exponent of the Brownian loop measure that
   arises in the Radon-Nikodym derivatives.
\item  $c \in (-\infty,1]$ is the {\em central charge}.  We will not use this parameter in this
paper, but $c$ comes up naturally in representations of
conformal field theories.
\item $d \in (1,3/2]$ is the Hausdorff dimension of the paths. This is also a parameter that will not be discussed in this paper.
\end{itemize}
We will tend to use $b$ as our parameter, but any of the quantities above
could equally well be used.

We have defined the measure $Q_{D,b,n}(\z,\w)$ for fixed $\z,\w$.  However,
if $D$ has a piecewise analytic boundary, we can define
\[   Q_{D,b,n} = \int_{(\p D)^n \times (\p D)^n} Q_{D,b,n}(\z,\w)
  \; d \z \, d\w. \]
In the case $b=1$, the scaling rule   implies that
$ Q_{D,1,n}$ is conformally invariant. 
               

\section{Preliminaries}

\subsection{$\SLE_\kappa$}\label{SLEdefnsect}

We review the definition
of  chordal $\SLE_\kappa$ in simply connected domains.
We will only consider the case $\kappa \leq 4$ in this paper, and let
$a = 2/\kappa$ as in~(\ref{paramdefns}).  Then $\SLE_\kappa$
in $\Half$ from $0$ to $\infty$ is the random curve
$\gamma:(0,\infty) \rightarrow \Half$ with $\gamma(0+) = 0$ satisfying
the following.  Let $g_t: \Half \setminus \gamma(0,t] \rightarrow
\Half$ be the unique conformal transformation with $g_t(z) -z
= o(1)$ as $z \rightarrow \infty$.  Then $g_t$ satisfies the
chordal Loewner differential equation
\begin{equation}  \label{loewner}
              \dot g_t(z) = \frac{a}{g_t(z) - U_t} , \;\;\;\;
   g_0(z) = z , 
\end{equation}
where the driving function $U_t$ is a standard one-dimensional Brownian motion. The
curve has been parametrized so that $\gamma(0,t]$ has
half-plane capacity $at$, i.e., so that
 $g_t$ has the expansion
\begin{equation}  \label{parameter}
   g_t(z) = z + \frac{at}{z} + O(|z|^{-2}), \;\;\;\;\; z \rightarrow
  \infty . 
\end{equation}

If $D$ is a simply connected domain and $z$, $w$ are distinct boundary
points of $D$, then chordal $\SLE_\kappa$ from $z$ to $w$ is defined
to be the conformal image of $\SLE_\kappa$ in $\Half$ under a conformal
transformation taking $0$ to $z$ and $\infty$ to $w$.  We consider
this as a probability measure on curves modulo reparametrization, and
denote this probability measure by $\mu^\#_{D,b}(z,w)$.

There are a number of special values for $\kappa$.  In this paper
we will look at special properties of $\kappa=2$. 
As we now show, the stochastic process 
$\dot g_t(z)$ 
is a martingale iff $\kappa=2$ (equivalently, $a=1$).

\begin{proposition}\label{MGprop}
Suppose that $g_t(z)$ satisfies the chordal
Loewner equation~(\ref{loewner}). If
$M_t=M_t^{z,a} = \dot g_t(z)$, then $M_t$ is a martingale if and only if $a=1$.
\end{proposition}

\begin{proof}
If $Z_t = g_t(z) -U_t$, then a simple calculation using It\^o's formula shows that
\[dZ_t = \frac{a}{Z_t}dt - dU_t,\]
so that $Z_t$ is a Bessel process with $d\langle Z \rangle_t = dt$.
Using It\^o's formula again therefore gives
\[dM_t = d \left(\frac{a}{Z_t}\right) = 
-\frac{a}{Z_t^2}dZ_t + \frac{a}{Z_t^3}d\langle Z \rangle_t = \frac{a}{Z_t^2}dU_t +
 \frac{a(1-a)}{Z_t^3}dt.\]
Hence, $M_t$ is a martingale iff $a=1$, in which case $dM_t = M_t^2 dU_t$. 
\end{proof}

\subsection{Brownian loop measure}   \label{Bloopsec}

We  review the Brownian loop measure as introduced
in~\cite{loopsoup}; see also~\cite[Section~5.7]{Lbook}.  (That paper 
considered the loop soup, which is
a Poissonian realization of the loop measure.  This loop soup is
very useful heuristically but is not needed for the definition
of the measure.)  A {\em rooted loop} of
time duration $t$ is a continuous function $\eta:[0,t] \rightarrow
\C$ with $\eta(0) = \eta(t)$.  It can also be considered as an
ordered pair $(t,\eta)$ where  $\eta:(-\infty,\infty) \rightarrow \C$ 
is a continuous function with period $t$.
An {\em unrooted loop} is an equivalence class of rooted loops under
the equivalence $(t,\eta) \sim (t,\eta_s)$ for all $s$
where $\eta_s(r) = \eta(s+r)$.  Each rooted loop can be written
uniquely as $(t,z,\bar \eta)$ where $t > 0$ is the time
duration, $z$ is the root, and $\bar \eta$ is a loop of time
duration $1$ rooted at the origin. More specifically,
$\eta$ is given by translation and Brownian scaling, 
$\eta(s) = z +  t^{1/2} \, \bar \eta(s/t)$.
  The {\em rooted Brownian loop measure}  
is given by choosing $(t,z,\bar \eta)$ from the measure
\[      \frac{1}{2 \pi t}  \, dt \times \mbox{ area } \times \nu_{BB}. \]
where  $\nu_{BB}$ denotes the probability measure associated to
Brownian bridges.  The {\em Brownian loop measure} is the measure
on unrooted loops obtained from this measure by ``forgetting the loops.''
The loop measure in a domain $D$ is exactly this measure
restricted to loops that stay in the domain. An important fact
is that the Brownian loop measure is conformally invariant: if
$f : D \rightarrow D'$ is a conformal transformation, then the image
of the loop measure on $D$ by $f$ is the loop measure on $D'$.
(If we are worrying about the parametrization of the loops, we have
to change the parametrization by the usual rule for conformal 
invariance of Brownian motion.  In this paper, the parametrization
will not be important, only the traces of the loops.)  The Brownian
loop measure can also be obtained as the scaling limit of a
simple measure on random walk loops; see Section~\ref{lerwsec}.

If $D \subsetneq \C$ is simply connected and $V_1,V_2$ are closed
sets,  
we let $m(D; V_1,V_2) < \infty$ denote the Brownian loop measure
of loops in $D$ that intersect both $V_1$ and $V_2$.  If $f:D \rightarrow
D'$ is a conformal transformation, then $m(f(D); f(V_1),f(V_2))
= m(D;V_1,V_2)$. The two main properties of the Brownian loop measure that
 we will use are conformal invariance, and the fact that if $D$ is simply connected, 
then $0 < m(D;V_1,V_2) < \infty$, provided that $V_1$, $V_2$ are disjoint 
closed sets with at least one $V_i$ compact; $\dist(V_1, V_2)> 0$;  and 
non-polar, i.e., such that with positive probability
a Brownian motion starting at $z \in D$ hits both of these sets.

\subsection{Boundary perturbation of $\SLE$}  \label{SLEperturb}

Suppose $D\subset \Half
 $ is a simply connected domain such that $\Half \setminus
D$ is bounded,  and $x \in \R$ with $\dist(x, \Half \setminus
D) > 0$.   For any such domain there
is a unique conformal transformation $\Phi = \Phi_{D,x}:
D \rightarrow \Half$ such that
$\Phi(x) =  0$, $\Phi(\infty) = \infty$, $\Phi'(\infty) = 1$.  Moreover,
Schwarz reflection tells us how to extend
$\Phi$  analytically in a neighborhood of $x$.  Schramm~\cite{schramm}
defined $\SLE_\kappa$  (modulo time parametrization) in $D$ from $x$ to $\infty$ to be 
the image of $\SLE_\kappa$ in $\Half$ from $0$ to $\infty$ under the
conformal map $\Phi^{-1}$.  By using properties of conformal
mappings and It\^o's formula (see, e.g.,~\cite[Section~4.6]{Lbook})
one can give an equivalent formulation (see, e.g.,~\cite[Section 4.3]{laplacian}). 
 Let $\gamma$ and $g_t$ be
as above, and let $D_t = g_t(D)$, $\Phi_t = \Phi_{D_t,U_t}$.
Consider the solution
to~(\ref{loewner}) where the driving function $U_t$ satisfies the
SDE
\begin{equation}  \label{nov19.1}
    dU_t =  b \, \frac{\Phi_t''(U_t)}{\Phi_t'(U_t)} \, dt
         +  dB_t, 
\end{equation}
where
$B_t$ is a standard one-dimensional Brownian motion. 
Then the distribution of the path $\gamma$ is exactly
that of  $\SLE_\kappa$ in
$D$.  In this formulation  we have parametrized the paths so that
$g_t$ satisfies~(\ref{parameter}).   This is not the same
parametrization as one would get if one started with an $\SLE_\kappa$
in $\Half$ from $0$ to $\infty$, say $\tilde \gamma$, and let
$\gamma(t) = \Phi^{-1} \tilde \gamma(t)$.  However, we only
define $\SLE_\kappa$ in domains up to time reparametrization, so this
is not an issue.  In fact, we could have used~(\ref{nov19.1}) as the
{\em definition} of $\SLE_\kappa$ in $D$, and then conformal invariance
of the process would be a {\em result}.

There is a another equivalent way of defining $\SLE_\kappa$   in $D$
 (for $\kappa \leq 4$)
that comes from a martingale introduced in~\cite{restriction}.  Let
$\gamma$ be an $\SLE_\kappa$ in $\Half$ from $x$ to $\infty$ parametrized
as in~(\ref{parameter}).  Let $T = T_D =
\inf\{t > 0: \gamma(t) \in \Half \setminus D\}$.  If $t < T$, then
we can let $D_t$ and $\Phi_t$ be as in the previous paragraph.  Define
$Y_t$ by
\[   Y_t = 1\{\gamma[0,t] \cap (\Half \setminus D) = \emptyset \}
         \, \exp\left\{\frac{a\lambda}{12} \int_0^t S \Phi_s(U_s) \, ds
\right\}, \]
where, as before, 
\[           \lambda = \lambda(\kappa) = \frac{(3a-1)(4a-3)}{2a} = 
\frac{(6-\kappa)(8-3\kappa)}{4\kappa}, \]
and $S$ denotes the Schwarzian derivative.  
In fact (see~\cite{restriction} and~\cite{loopsoup}), 
\[                      - \frac{a}{12} \int_0^t S \Phi_s(U_s) \, ds
 =m(\Half; \Half \setminus D , \gamma(0,t]).\]
(The factor $a$ comes from the fact that $\gamma$ is parametrized as
 in~(\ref{parameter}).)
Hence, we can write
\begin{equation}  \label{nov21.2}
     Y_t = 1\{\gamma[0,t] \cap (\Half \setminus D) = \emptyset \}
          \, \exp\left\{-\lambda \, m(\Half; \Half \setminus D , \gamma(0,t])
  \right\}. 
\end{equation}
In~\cite{restriction} it was shown that  $M_t = \Phi_t'(U_t)^b \,
Y_t$ is a local martingale.  It is easy to check that $
\Phi_t'(U_t) \leq 1$ which implies
that  $M_t$ is a bounded martingale for $\lambda \geq 0$
 $(\kappa \leq 8/3)$.  In this case, 
Girsanov's theorem can be
used to find the effective SDE for paths weighted by the random
variable $M_\infty$.   The resulting equation  is~(\ref{nov19.1}).
Since this is the same as $\SLE_\kappa$ in $D$, and we know that (for
$\kappa \leq 4$) $\SLE_\kappa$ in $D$ never leaves $D$, we conclude 
that with respect to the weighted measure, $T = \infty$.
What this says is that for $\kappa \leq 8/3$, the probability measure on 
paths given by starting with $\SLE_\kappa$ from $x$ to $\infty$ in $\Half$,
and then taking the Radon-Nikodym derivative $M_\infty/M_0$, is the same
as $\SLE_\kappa$ in $D$.  Since $M_0 = \Phi'(x)^b$ and $M_\infty = Y_\infty$,
 we conclude that
$\E[Y_\infty] = \Phi'(x)^b$.  

For $8/3 < \kappa \leq 4$, we have $\lambda < 0$ and hence $Y_t$ is
not bounded.  However, the same conclusions hold.  The reason is as follows.
For a fixed $\gamma:(0,\infty) \rightarrow \Half$ parametrized as
in~(\ref{parameter}) with $\gamma(t) \rightarrow \infty$ and $\dist(\gamma(0,\infty),
\Half \setminus D)>0$, standard conformal mapping estimates show that
$Y_t \rightarrow Y_\infty < \infty$.  In other words, the only problem
that one can have is that $Y_t$ blows up as $\gamma(t) \rightarrow \Half \setminus D$.
If we consider a stopping time of the form $\tau_n = \inf\{t: Y_t \geq n\}$,
then we know that $Y_{t \wedge \tau_n}$ is a martingale.  If we use
Girsanov to weight the paths by this martingale, the paths satisfy~(\ref{nov19.1}).  
We know (for $\kappa \leq 4$)
that paths that are given by this SDE never
reach $\Half \setminus D$ and go to infinity.  Hence, $M_t$ is actually a uniformly
integrable martingale with $\E[M_\infty] = \E[M_0]
= \Phi'(x)^{b} \leq 1$. 

As a prelude to the next section, we can see that we could define the
$\SLE_\kappa$ measure in $D$ to be the measure on paths whose
Radon-Nikodym derivative with respect to $\SLE_\kappa$ in $D$ is
 $Y_\infty$.  This is a finite measure with total measure
$\Phi'(x)^b$.  If we normalize to make this a probability measure,
then the corresponding measure is that of $\SLE_\kappa$ in $D$ from
$x$ to $\infty$.


\section{$\SLE$ measure}

\subsection{One path ($n=1$)}

The configurational approach
is  to view chordal $\SLE_\kappa$ as not just a probability
measure on paths connecting two specified points on the boundary, but
rather as a finite
measure on paths that 
when normalized gives  
 chordal $\SLE_\kappa$ as defined by Schramm.
Recall that
we  restrict our consideration to $\kappa \leq 4$.
Suppose $D$ is a simply connected domain and
$z$, $w$ are distinct boundary points at which $\p D$ is 
analytic.  Then the $\SLE_\kappa$ measure on   
paths in $D$ from $z$ to $w$ is defined to be
\[   Q_{D,b}(z,w) =  H_{D,b}(z,w) \; \mu^\#_{D,b}(z,w) \]
where $\mu^\#_{D,b}(z,w)$ is the $\SLE_\kappa$
probability measure, and $H_{D,b}(z,w) =|Q_{D,b}(z,w)|$ 
is determined by the scaling rule~(\ref{hinvariance}) from
the kernel
\[              H_b(x,y) = H_{\Half,b}(x,y) ={|y-x|^{-2b}}  \]
as in~(\ref{3}). The conformal covariance rule
\[  f \circ Q_{D,b}(z,w)  =  |f'(z)|^b \, |f'(w)|^b
 \; Q_{f(D),b}(f(z),f(w)) \]
follows immediately from the conformal invariance
of  $\mu^\#_{D,b}(z,w)$ and the scaling rule~(\ref{hinvariance})
for $ H_{D,b}(x,y) $.
  An important property of the measure
is the boundary perturbation rule.

\begin{proposition}[Boundary Perturbation] \label{nov22.prop1}
Suppose $D \subset D' \subsetneq \C$ are simply
connected domains.  Suppose that $z$, $w$ are distinct points of
$\p D$; $\p D$ is analytic in neighborhoods of $z$, $w$; and
$\p D$, $\p D'$ agree in neighborhoods of $z$, $w$.  Then
$Q_{D,b}(z,w)$ is absolutely continuous with respect to
$Q_{D',b}(z,w)$ with Radon-Nikodym derivative
\[    Y_{D,D',b}(z,w)(\gamma)   =    1\{ \gamma \subset D\}
  \,\exp\big\{-\lambda \, m(D';\gamma,D' \setminus D)\big\} . \] 
   In particular,
the Radon-Nikodym derivative is a conformal invariant,
and moreover
\begin{equation}  \label{nov21.3}
   H_{D,b}(z,w) \leq H_{D',b}(z,w). 
\end{equation}
\end{proposition}

\begin{proof}
It suffices to prove this for $D' = \Half$ and this
is discussed in Section~\ref{SLEperturb}. 
\end{proof}

\begin{remark} If $f : D' \rightarrow f(D')$
is a conformal transformation, then \[ Y_{f(D),f(D'),b}(f(z),f(w))
(f \circ \gamma) =Y_{D,D',b}(z,w)(\gamma).\]
It follows that the Radon-Nikodym
derivative
\[  Y_{D,D',b}(z,w) = \frac{d Q_{D,b}(z,w)}
         {dQ_{D',b}(z,w)} \]
makes sense even for non-smooth boundaries provided that $D 
\subset D'$ and the boundaries agree near $z$ and $w$.  If $f:
D' \rightarrow \Half$ is a conformal transformation with $f(z) = 0,
f(w) = \infty$, then $f(D' \setminus D)$ is a bounded subset
of $\Half$ bounded away from $0$. 
\end{remark}

\subsection{Multiple paths}\label{multpathsect}

We will now define the measures  $Q_{D,b,n} $ for positive
integers $n$.   Let $\z = (z_1,\ldots,z_n)$, $\w = (w_1,\ldots,
w_n)$ denote $n$-tuples of distinct points in $\p D$ ordered counterclockwise as
in the introduction.  Then
$Q_{D,b,n}(\z,\w)$,
 the $n$-path $\SLE_\kappa$ measure in $D$,
 is defined to be the measure
that  is  absolutely continuous with respect to
the product measure 
\[Q_{D,b }(z_1,w_1) \times \cdots\times Q_{D,b }(z_n,w_n)\]   
with Radon-Nikodym derivative 
$  Y({\bar \gamma}) =  Y_{D,b,\z,\w} (\gamma^1,\ldots,\gamma^n)$ given by
\begin{equation}  \label{nov21.1}
Y({\bar \gamma}) = 1\{ \gamma^k \cap \gamma^l
= \emptyset, \; 1 \leq k < l \leq n \} \;
\exp\left\{-\lambda \sum_{k=1}^{n-1}
             m(D;\gamma^k,\gamma^{k+1})\right\}. 
\end{equation}
If $\lambda \geq 0$, it is immediate that
$Q_{D,b,n}(\z,\w) $ is a finite measure.  As we will see after
we prove the next proposition, 
this is also true for $\lambda \geq -1/2$;
see~(\ref{nov21.4}). 

Note that $Y$ is conformally invariant.  Consequently,
the measures  $Q_{D,b,n} $ satisfy the conformal
covariance relation~(\ref{confcov}). 
  Other important properties are
given in the next two propositions.

\begin{proposition}[Cascade Relation] For $1\leq j \leq n$, if
  \[ \z = (z_1,\ldots,z_n), \;\;\;\; \w = (w_1,\ldots,
w_n), \;\;\;\; \hatgamma = (\gamma^1,\ldots,\gamma^{j-1},
\gamma^{j+1},\ldots,\gamma^n),\] \[ \hat \z = (z_1,\ldots,z_{j-1}, z_{j+1},\ldots,
z_n) ,\;\;\;\; \hat \w = (w_1,\ldots,w_{j-1},w_{j+1},\ldots,w_n), \]
then the marginal measure on $\hatgamma$  in $Q_{D,b,n}(\z,\w)$
is absolutely continuous with respect to $Q_{D,b,n-1}(\hat \z,
\hat \w)$ with Radon-Nikodym derivative
$H_{\hat D,b}(z_j,w_j)$. 
Here $\hat D$ is the subdomain of $D \setminus \hatgamma$ whose boundary includes
$z_j$, $w_j$.  (If $z_j$, $w_j$ are contained in the boundary of
different components, then $H_{\hat D,b}(z_j,w_j) = 0$.)
Moreover, the conditional distribution of $\gamma^j$ given
$\hatgamma$
is that of $\SLE_\kappa$ from $z_j$ to $w_j$ in $\hat D$.
\end{proposition}

\begin{proof}
 Let $\hat Y(\hatgamma)=Y_{D,b,\hat \z, \hat \w}(\hatgamma)$ denote 
the density of $Q_{D,b,n-1}(\hat \z,
\hat \w)$ with respect to the appropriate product measure
as defined in~(\ref{nov21.1}).
Then 
\[   Y 
 =  \hat Y \, 1\{ \gamma^l \cap \gamma^k
= \emptyset,\; 1 \leq l < k \leq n \} \;
\exp\left\{- \lambda \hat m\right\} 
             \] where
\[ \hat m = m(D; \gamma^{j-1}, \gamma^j) +
m(D;\gamma^j,\gamma^{j+1}) -
  m(D; \gamma^{j-1}, \gamma^{j+1}) \]
Simple inclusion-exclusion shows that
$\hat m$ is the measure of loops in $D$ that intersect both
$\gamma^j$ and $D \setminus \hat D$.  The result now
follows from~(\ref{nov21.1}) and the discussion in 
Section~\ref{SLEperturb}. 
\end{proof}

\begin{remark}
Let $H_{D,b,n}(\z,\w) = |  Q_{D,b,n}(\z,\w)|$. It follows from
the proposition and~(\ref{nov21.3}), that
\[  H_{D,b,n}(\z,\w) \leq  H_{D,b,n-1}(\hat \z,\hat \w)
  \, H_{D,b}(z_j,w_j). \]
By iterating, we get 
\begin{equation}  \label{nov21.4}
   H_{D,b,n}(\z,\w) \leq \prod_{j=1}^n  H_{D,b}(z_j,w_j), 
\end{equation}
so that we conclude $Q_{D,b,n}(\z,\w)$ is a finite measure for all $b \ge 1/4$, or equivalently, all $\lambda \ge -1/2$.
\end{remark}

\begin{proposition}[Boundary Perturbation]
Suppose $D \subset D' \subsetneq \C$ are simply
connected domains.    Then
$Q_{D,b,n}(\z,\w)$ is absolutely continuous with respect to
$Q_{D',b,n}(\z,\w)$ with Radon-Nikodym derivative
\[  Y_{D,D',b,n}(\z,\w)(\bar \gamma)  = 
    1\{ \gamma^j \subset D, \;
j=1,\ldots,n\} \; 
\exp\{-\lambda \,m(D';\gamma^1 \cup \cdots \cup
  \gamma^n, D' \setminus D)\}. \]
  In particular,
the Radon-Nikodym derivative is a conformal invariant.
\end{proposition}

\begin{proof}
It follows from Proposition~\ref{nov22.prop1} that 
the Radon-Nikodym derivative of
\[Q_{D,b }(z_1,w_1) \times \cdots\times Q_{D,b }(z_n,w_n)\]   
with respect to
\[Q_{D',b }(z_1,w_1) \times \cdots\times Q_{D',b }(z_n,w_n)\]   
is  
 \[
     1\{ \gamma^j \subset D, \;
j=1,\ldots,n\} \;
\exp\left\{  - \lambda
\sum_{k=1}^n m(D';\gamma^k,D' \setminus D ) \right\}. \] 
Hence by~(\ref{nov21.1}),  $ Y_{D,D',b,n}(\z,\w)(\bar \gamma)$
equals
 this quantity times
\[    
\exp\left\{-\lambda \sum_{k=1}^{n-1}
             [m(D;\gamma^k,\gamma^{k+1}) - m(D';
   \gamma^k,\gamma^{k+1})]\right\} . \]
But $m(D';\gamma^k,\gamma^{k+1}) - m(D;
   \gamma^k,\gamma^{k+1})$ denotes the measure of the set
of loops in $D'$ that intersect all three
of $\gamma^k, \gamma^{k+1}$, and $D' \setminus D$.  Therefore,
for $k=1,\ldots, n-1$,
\[ m(D';\gamma^k,D' \setminus D ) - [m(D';\gamma^k,\gamma^{k+1}) - m(D;
   \gamma^k,\gamma^{k+1})] \]
denotes the measure of the set of loops in $D$
 that intersect  $\gamma^k$ and $D' \setminus D$ but do not 
intersect $\gamma^{k+1}$.  For topological reasons, they cannot
possibly intersect $\gamma^l$ for $l \geq k+1$.  For any loop in
$D'$ that intersects both $\gamma^1 \cup \cdots \cup \gamma^n$, there is
a largest $k$ such that the loop intersects $\gamma^k$ but not
$\gamma^{k+1},\ldots,\gamma^n$.  This is accounted for in the $k$th term
of the sum and in no other.
\end{proof}

Recall that 
  $H_{D,b,n}$ satisfies the  scaling rule
\[   H_{D,b,n} (\z,\w) = |f'(\z)|^b \, |f'(\w)|^b \,
     H_{f(D),b,n}(f(\z),f(\w)) . \]
It also satisfies the following cascade relation.  Suppose
$\z = (z_1,\hat \z)$, $\w = (w_1,\hat \w)$ and all the points
of $\hat \z$, $\hat \w$ lie on the same arc connecting $z_1$
and $w_1$ on $\p D$.  Then,
\[  H_{D,b,n}(\z, \w) =  Q_{D,b}(z_1,w_1) \, [\,
       H_{D(\gamma),b,n-1}(\hat \z, \hat \w)\,] \]
where $D(\gamma)$ is the connected component of $D \setminus \gamma$
whose boundary contains  $\hat \z$, $\hat \w$ and we use the
shorthand notation $\mu \,[\,f\,] = \int f \, d\mu$.
If $f:D \rightarrow \Half$ is a conformal transformation with
$f(z_1) = 0$, $f(w_1) = \infty$, and $\hat \z$, $\hat \w$ are sent to
the positive real axis, then
\begin{equation}\label{dec16eq1}
   H_{D,b,n}(\z,\w) =  |Q_{D,b}(z_1, w_1)|
 \, |f'(\hat \z)|^b \, |f'(\hat \w)|^b \,
             \E[ H_{(\Half \setminus \gamma)_+,b,n-1}(f(\hat \z), f(\hat \w))],
\end{equation}
where in this case the expectation is over $\SLE_\kappa$ from
$0$ to $\infty$ in $\Half$ and $(\Half \setminus \gamma)_+$ denotes the
component of $\Half \setminus \gamma$ to the right of $\gamma$.
This shows that $Q_{D,b,n}$ for all $n$ can be computed, at least
in principle, from computations of $\SLE_\kappa$ in $\Half$.

We will also define the \emph{partition function}
\begin{equation}
\label{nov23.1}
  \tilde H_{D,b, n}(\z,\w) = \frac{H_{D,b ,n} (\z,\w)}
     {H_{D,b}(z_1,w_1) \cdots H_{D,b}(z_n,w_n)}. 
\end{equation}
Note that this is a conformal {\em invariant},
i.e., \[  \tilde H_{D,b,n} (\z,\w) = \tilde
H_{f(D),b,n}(f(\z),f(\w)).\]
In particular $H_{D,b,n}$,and $\tilde H_{D,b,n}$ are determined by the values of
\[  H_{b,n}(\x,\y) = H_{\Half,b,n}(\x,\y) , \]
for $\x$, $\y \in \R^n$. We also write $\tilde H_{b,n}$ for $\tilde H_{\Half,b,n}$.

\begin{remark}
One case that is known is 
$\SLE_2$, corresponding to \emph{loop-erased walk} ($a=b=1$), in which case
\[   H_{1,n}(\x,\y) = \det[(y_k - x_j)^{-2}]_{1 \le j,k \le n} . \]
This is a consequence of Fomin's identity and the $n=1$ case.
We describe this in more detail in Section~\ref{fominscalesect}.
\end{remark}

If $0 < x_1 < \cdots < x_n < y_n < \cdots <y_1 < \infty$, we
let
\[\tilde H_{b, n}^*(\x,\y) 
 =    \tilde H_{b, n+1}((0,\x), (\infty,\y))
   =  \lim_{w \rightarrow \infty}
                  \tilde H_{b, n+1}((0,\x), (w,\y)),\]
\[  H_{b, n}^*(\x,\y)  = \left[\prod_{j=1}^n H_{b,1}(x_j,y_j
  ) \right] \, \tilde H_{b, n}^*(\x,\y) 
     =    \lim_{w \rightarrow \infty}
     w^{2b} \,  H_{b, n+1}((0,\x), (w,\y)). \]
These functions satisfy the following scaling rules for $r> 0$:
\begin{equation}\label{scalerule}
    \tilde H_{b,  n}^*(r\x,r\y)  = \tilde H_{b, n}^*( \x, \y) ,
\;\;\;\;  H_{b, n}^*(r\x,r\y) = r^{-2bn} H_{b,  n}^*( \x, \y).
\end{equation}
 
\subsection{The partition function for two paths} \label{partsec}

In this section we calculate $H_{b,1}^*(x,y) = (y-x)^{-2b} \,
  \tilde H_{b,1}^*(x,y) = (y-x)^{-2b} \tilde H_{b,2}((0,x),(\infty,y))$ 
for $0 < x < y <\infty$. The formula in the following proposition has appeared
 previously; see, for example,~\cite{dubedat}. As usual, we will write just $H^*(x,y)$, 
 $\tilde H^*(x,y)$ for $H_{b,1}^*(x,y)$,  $\tilde H_{b,1}^*(x,y)$, respectively. 
We will use this result in Section~\ref{directderivsect} to derive the scaling 
limit of Fomin's identity.

\begin{proposition} \label{dec1.prop1}
 If $b \geq 1/4$, then
\begin{equation}  \label{nov19.2}
H^*_{b,1}(x,y) = (y-x)^{-2b} \, \frac{\Gamma(2a) \, \Gamma(6a-1)}
            {\Gamma(4a) \, \Gamma(4a-1)} \, (x/y)^a \,
                F(2a,1-2a,4a; x/y) =\frac{\phi_b(x/y)}{(y-x)^{2b}}, 
\end{equation}
where $F$ denotes the hypergeometric function
and $a = (2b+1)/3$.
\end{proposition}

\begin{proof}
By the scaling rule~(\ref{scalerule}), we know that $\tilde H^*(x,y)
= \phi(x/y)$ for some function $\phi = \phi_b$ of one variable. 
We will parametrize $\SLE_{\kappa}$ so that it satisfies~(\ref{loewner})
where the driving function $U_t = - B_t$ is a standard one-dimensional Brownian motion.   
Let $0 < x < y$, and let $X_t = g_t(x) + B_t$, $Y_t = g_t(y) + B_t$ so that
\[  dX_t = \frac{a}{X_t} \, dt + dB_t, \;\;\;\;\;
      dY_t = \frac{a}{Y_t}\, dt +dB_t \]
as in the proof of Proposition~\ref{MGprop}. 
Let $(\Half\setminus \gamma)_+$ denote the component of $\Half \setminus \gamma$ whose
boundary includes $x$, $y$.  If $g:(\Half\setminus \gamma)_+ \rightarrow \Half $ is a
conformal transformation with $g(\infty) = \infty$, then
\[  H_{(\Half\setminus \gamma)_+,b}(x,y) = \frac{g'(x)^b \, g'(y)^b}{[g(y) -
    g(x)]^{2b}} \]
as a consequence of~(\ref{hinvariance}) and~(\ref{3}).
Using this we see that if 
\[             J_t = \frac{g_t'(x)^b \, g_t'(y)^b}{[g_t(y) - g_t(x)]^{2b}} 
\;\;\; \text{ and } \;\;\;  J_\infty = \lim_{t \rightarrow \infty} J_t \]
then from~(\ref{dec16eq1}) and~(\ref{nov23.1})
\begin{equation}\label{nov17.01}         
\tilde H^*(x,y) =  (y-x)^{2b} \, \E^{x,y}[J_\infty]
\end{equation}
so that $H^*(x,y)=\E^{x,y}[J_\infty]$.
We will now give the differential equation for  $\tilde H^*(x,y)$.  Since
\[     \p_t \, [\log g_t'(x)] = -\frac{a}{X_t^2}, \;\;\;\;\;
              \p_t \,[\log g_t'(y)] = - \frac{a}{Y_t^2}, \;\;\;\; \text{ and }\]
\[      \p_t\,[\log(g_t(y) - g_t(x))]  =
                 \frac{1}{Y_t - X_t} \, \left[\frac{a}{Y_t} - \frac
   a{X_t} \right] = -\frac{a}{X_t \, Y_t},  \]
we see that
\begin{align*}
J_t &= \exp\left\{ \log J_t \right\}
=J_0 \, \exp\left\{\int_0^t \p_s[\log J_s] \; ds \right\}
=\frac{1}{(y-x)^{2b}} \, \exp\left\{-ab \int_0^t \left( \frac{1}{X_s} - \frac{1}{Y_s}
 \right)^2 \; ds \right\}.
\end{align*}
Hence,~(\ref{nov17.01}) gives
\[    \tilde H^*(x,y) = \E^{x,y} \left[
 \exp\left\{-ab \int_0^{\infty} \left( \frac{1}{X_s} - \frac
   1 {Y_s} \right)^2 \; ds \right\}\right].\]
It now follows from the (usual) Markov property that $J_t \, \tilde H^*(X_t,Y_t)$ is a 
martingale. That is, 
if $M_t = \E^{x,y}[J_{\infty}| \mathcal{F}_t]$ so that $M_t$ is a martingale, then
\begin{align*}
M_t 
&=\E^{x,y} \left[ \;\exp\left\{-ab \int_0^{\infty} \left( \frac{1}{X_s} - \frac{1}{Y_s}
 \right)^2 \; ds \right\} 
\; \bigg| \;\mathcal{F}_t\;\right]\\
&=\exp\left\{-ab \int_0^{t} \left( \frac{1}{X_s} - \frac{1}{Y_s} \right)^2 \; ds \right\} 
\E^{x,y} \left[ \; \exp\left\{-ab \int_{t}^{\infty} \left( \frac{1}{X_s} - \frac{1}{Y_s}
 \right)^2 \; ds \right\} 
\; \bigg| \;\mathcal{F}_t\;\right]\\
&=J_t\, \tilde H^*(X_t,Y_t).
\end{align*}
It\^o's formula now implies that 
\begin{equation}\label{ItoEqn}
  -ab\left(\frac 1x - \frac 1y \right)^2 \, \tilde H^*  + \frac ax \, 
    \p_x  \tilde H^*  
+ \frac ay \, \p_y   \tilde H^*  + \frac 12 \, \p_{xx}
     \tilde H^*     + \frac 12 \, \p_{yy}  \tilde H^* 
   + \p _{xy} \tilde H^* = 0 .
\end{equation}
Since $\tilde H^*(x,y) = \phi(x/y)$, we find
\[   \p_x \tilde H^* = y^{-1} \, \phi'(x/y), \;\;\;   \p_y \tilde H^* = - x \, y^{-2} \, \phi'(x/y), \;\;\;
   \p_{xx} \tilde H^* = y^{-2} \, \phi''(x/y), \]
\[   \p_{yy} \tilde H^* = 2\,  x \, y^{-3} \, \phi'(x/y) + x^2 \, y^{-4} \,   \phi''(x/y), \;\;\;
 \p_{xy} \tilde H^* = - y^{-2} \, \phi'(x/y) - x \, y^{-3} \, \phi''(x/y),\]
so that after substituting into~(\ref{ItoEqn}), multiplying by $y^2$, letting $u=x/y$,
combining terms, and recalling that $b = (3a-1)/2$, we have
\begin{equation}\label{hypergeomFomin}
 u^2 \, (1-u)^2 \,\phi''(u)  +  2 \,
 u \, (a-u + (1-a) \, u^2) \, \phi'(u) -a(3a-1) (1-u)^2 \, \phi(u) =0.
\end{equation}
Note that $0<u<1$, and let $\psi(u) = u^{-a} (1-u)^{1-4a} \phi(u)$ so that~(\ref{hypergeomFomin}) implies
\begin{equation}\label{hypergeomFomin3}
 u \, (1-u) \,\psi''(u)  +  (4a-8au) \, \psi'(u) -2a(6a-1)\, \psi(u) =0.
\end{equation}
We see that~(\ref{hypergeomFomin3}) is the well-known hypergeometric differential equation (see equation~(9.5.5) of~\cite{lebedev})
which has two linearly independent solution; hence
\[\psi(u) = C_1 F(2a,6a-1,4a;u) + C_2 u^{1-4a} F(1-2a,2a,2-4a;u)\]
where $C_1$ and $C_2$ are arbitrary constants.
Equation (9.5.3) of~\cite{lebedev} gives
$        F(2a,6a-1,4a;u) = (1-u)^{1-4a} \, F(2a,1-2a,4a;u) $
so that the general solution of~(\ref{hypergeomFomin}) is therefore
\[\phi(u)=    C_1  u^a \,  F(2a,1-2a,4a;u)
   +C_2 u^{1-3a} \, (1-u)^{4a-1} \, F(2a,1-2a,2-4a;u).\]
Since $\phi(u) \to 0$ as $u \to 0+$, and $\phi(u) \to 1$ as $u \to 1-$, we find $C_2=0$ and
\[ C_1^{-1}=           F(2a,1-2a,4a;1) = \lim_{u\to 1-}F(2a,1-2a,4a;u) =\frac{\Gamma(4a) \, \Gamma(4a-1)}
                {\Gamma(2a) \, \Gamma(6a-1)}, \]
so that
\begin{equation*}
 \phi(u) =\phi_b(u)= \frac{\Gamma(2a) \, \Gamma(6a-1)}
            {\Gamma(4a) \, \Gamma(4a-1)} \, u^a \,
                F(2a,1-2a,4a;u).\qedhere
\end{equation*}
\end{proof}

\begin{remark}
If $x > 0$ is fixed, then 
\begin{equation}  \label{dec24.2}
   \lim_{y \rightarrow \infty}
  \, y^{a+2b}  H_{b,1}^*(x,y)  = \frac{\Gamma(2a) \, \Gamma(6a-1)}
            {\Gamma(4a) \, \Gamma(4a-1)}  \, x^a. 
\end{equation}
\end{remark}

\begin{remark}
If $\tilde H_{b,2} = \tilde H_{\Half,b,2}$
is as defined in~(\ref{nov23.1}), then this proposition implies
for $x_1 < x_2 < y_2 < y_1$ that
\begin{equation}  \label{nov23.2}
\tilde H_{b,2}(x_1,x_2,y_1,y_2) =
   \phi_b(1-q), 
\end{equation} 
where
\[  q = q(x_1,x_2,y_2,y_1) = 
\frac{(y_1 - x_1) \, (y_2 - x_2)}{(y_1 - x_2) \,
   (y_2 - x_1)}, \]
denotes the cross-ratio.
As pointed out in~\cite{dubedat}, if the quantity on the left of~(\ref{nov23.2}) is
to be a conformal invariant, it should be a function of $q$.  If $x_1 = 0$,
$x_2 = x$, $y_2 = y$, $y_1 = \infty$, then $q = 1-(x/y)$ as expected.
\end{remark}

\begin{remark}
If $1-2a = -k$ where $k$ is a non-negative integer,
then $F(2a,1-2a,4a;  x/y)$ is a $k$th degree polynomial.  These
values correspond to
\[            b = \frac{3k+1}{4}. \]
In particular, $\phi_{1/4}(u)=\sqrt{u}$ and  $\phi_1(u)=2u-u^2=1-(1-u)^2$.
Several other special cases of~(\ref{nov19.2}) are listed in the following table.
\begin{center}
\renewcommand{\arraystretch}{2.2}
\begin{tabular}{|c|c|c|c|}\hline
$a=$  &$b=$  &$\kappa=$  &$H_{b,1}^*(x,y) = $ \\ \hline\hline
$2$   &$5/2$ &$1$        &$(y-x)^{-5}\;  (x/y)^{2}\, [6 - 9(x/y) +5(x/y)^2-(x/y)^3]$\\\hline
$3/2$ &$7/4$ &$4/3$      &$(y-x)^{-7/2}\;  (x/y)^{3/2}\, [(7/2) - (7/2)(x/y) +(x/y)^2]$ \\\hline
$1$   &$1$   &$2$        &$(y-x)^{-2} \; (x/y) \, [2 - (x/y)]$ \\\hline
$1/2$ &$1/4$ &$4$        &$(y-x)^{-1/2} \;  (x/y)^{1/2}$ \\\hline
\end{tabular}
\end{center} 
\end{remark}

\begin{remark}
Our argument is essentially the same as the calculation of \emph{Cardy's formula} for $\SLE_\kappa$.
An examination of the proof will show that the hypergeometric solution is valid for $a>1/4$, or, equivalently, $\kappa < 8$.
The only reason we restrict to $b \ge 1/4$ in the proposition is to have $\kappa \le 4$. In fact, if $\kappa=6$ (corresponding to $b=0$, or $a=1/3$) then we recover Cardy's formula, namely
\[ \phi_0(u) = \frac{\Gamma(2/3)}{\Gamma(4/3) \, \Gamma(1/3)} \, u^{1/3} \,F(2/3, 1/3, 4/3; u). \]
A proof for the $b=0$ case only may be found in either~\cite{Car} or~\cite[Proposition~6.33]{Lbook}.
\end{remark}


\section{Multiple $\SLE$s}

As was shown in Section~\ref{multpathsect}, the definition of $Q_{D,b,n}(\z,\w)$ is fairly straightforward. Since $H_{D,b,n}(\z,\w) = | Q_{D,b,n}(\z,\w)|$, we can write
\[   Q_{D,b,n} (\z,\w) = H_{D,b,n}(\z,\w) \,
 \mu^\#_{D,b,n}(\z,\w), \]
for some probability measure $\mu^\#_{D,b,n}(\z,\w)$
 on $n$-tuples of
paths that can be called $\SLE_\kappa$ from
$\z$ to $\w$ in $D$.  This measure will be absolutely
continuous with respect to the product measure given
by taking $n$ independent $\SLE_\kappa$ paths 
$\gamma^1,\ldots,\gamma^n$ where $\gamma^j$ goes from
$z_j$ to $w_j$.  In the present section, we use Girsanov's theorem to give
the SDE.  

Suppose $G(t,x_0, x_1,\ldots,x_n)$ is a smooth,
non-negative function
(at least on the set $\{x_j \neq x_k: 0 \leq j < k
\leq n \}$), 
 and suppose $\gamma:(0,\infty) \rightarrow
\Half$ is an $\SLE_\kappa$ curve with corresponding
maps $g_t$ parametrized as in~(\ref{loewner}).  We allow the function $G(t,x_0, x_1,\ldots,x_n)$ to be
a random function adapted to the Brownian filtration,
but we assume that with probability one $G$ is 
$C^1$ in $t$.
Let $X^j_t = g_t(x_j) $ where
$W_t$ is a standard Brownian motion, and
note that 
\[    \dot X^j_t = \frac{a}{X^j_t-W_t}  . \]
If $J_t = G(t,W_t,X_t^1,\ldots,X_t^n)$, then It\^o's formula gives
\[  dJ_t = J_t \, \left[ {R_t}   \; dt
    +\frac{\p_0 G(t,W_t, \bar X_t)}{G(t,W_t, \bar X_t)}  \, dW_t\right], \]
where
\[ R_t = G(t,W_t, \bar X_t)^{-1} \, \left[ \dot G(t,W_t,\bar X_t) + 
\frac 12 \,  \p_{00}G  (t,W_t,\bar X_t) +
a\sum_{j=1}^n
        \frac{  \p_jG(t,W_t,  \bar X_t)}
              {X_t^j - W_t} \right]\]
with $\bar X_t = (X_t^1,\ldots, X_t^n)$ and $\p_j = \p_{x_j}$, $\p_{00}=\p_{x_0x_0}$.
If
\[   M_t = \exp\left\{-\int_0^t R_s \, ds  \right\} \, J_t,\]
then $M_t$ is a local martingale satisfying
\[   dM_t =  M_t\, \frac{\p_0G  (t,W_t, \bar X_t)}
   {G(t,W_t, \bar X_t)}  \, dW_t .\]
If $M_t$ is a martingale, Girsanov's
theorem says that if we weight the paths by $M_t$
then this is the same as if $W_t$ satisfied
\begin{equation}  \label{nov18.1}
   dW_t =  \frac{\p_0G(t,W_t, \bar X_t)}
   {G (t,W_t, \bar X_t)} \, dt + dB_t, 
\end{equation}
for a standard Brownian motion $B_t$.  Hence
an $\SLE_\kappa$ weighted by $G$ is the same as
a solution to the chordal Loewner equation~(\ref{loewner})
where the driving function $W_t$ satisfies~(\ref{nov18.1}).
Even if $M_t$ is not a martingale, we can get a similar
expression for a stopped process $M_{t \wedge \tau}$ if
the stopped process is a continuous martingale.

\subsection{$\SLE$ in a domain $D$}

 Suppose that $D \subset \Half$ is a simply
connected domain with $\Half \setminus D$ bounded and
bounded away from zero.  Let $D_t = g_t(D)$ and let
$\Phi_t$ be a conformal transformation of $D_t$ onto $\Half$
with $\Phi_t(\infty) = \infty$, $\Phi_t'(\infty) = 1$.
Let $G(t,x) = \Phi_t'(x)^b$.  Then if we weight the paths
of $\SLE_\kappa$ by this function, the paths satisfy
\begin{equation}  \label{dec1.1}
  dW_t = b \, \frac{\Phi_t''(W_t)}{\Phi_t'(W_t)} \, dt
    + d B_t. 
\end{equation}
This is the same equation as given in Section~\ref{SLEperturb}
 for $\SLE_\kappa$ in $D$
when one takes the conformal image of $\SLE_\kappa$ from
$\Half$ to $D$.  In other words, another way of viewing $\SLE_\kappa$
in $D$ is as $\SLE_\kappa$ weighted by $\Phi_t(W_t)^b$.

A similar example is $\SLE_\kappa$ from $0$ to $y$ where $0 < y < \infty$.
This measure on paths is absolutely continuous with respect 
to $\SLE_\kappa$ from $0$ to $\infty$ provided that the paths are stopped at a stopping
time smaller than the first visit to $y$.   At that point the paths
become singular.  (Although this is somewhat imprecise, it is not too difficult
to make it precise.)  The equation~(\ref{dec1.1}) can be applied
to see how the process evolves.  If $t$ is smaller than the hitting
time of $y$, we define $g_t$ as before and let $X_t = g_t(y)$.  Then,
$\Phi_t(z) = (z-X_t)^{-1}$ is a conformal map taking $g_t(\Half
\setminus \gamma(0,t])$ (= $\Half$ in this case) to $\Half$, mapping
$X_t$ to $\infty$.  Using~(\ref{dec1.1}), we get
\begin{equation}  \label{dec1.2}
  dW_t =  \frac{2b}{X_t - W_t} \, dt  + d B_t, \;\;\;\;\;\;
   dX_t = \frac{a}{X_t - W_t} \, dt.
\end{equation}
This is an example of an $\SLE(\kappa,\rho)$ process as introduced
in~\cite{restriction}.
Note that
\[   d(W_t - X_t) =  \frac{(1-4b)/3}{W_t - X_t} \, dt + dB_t , \]
so comparison with a Bessel process shows that the path
has finite lifetime.  (Recall that we are parametrizing the curve
by its half-plane capacity --- if $\gamma$ is a simple curve from
$0$ to $y$, then the path has finite capacity.)  We can derive the
same process by taking $\SLE_\kappa$ from $0$ to $\infty$ and weighting
the paths by $G(t,W_t,X_t) =
H_b(W_t,X_t) = (X_t - W_t)^{-2b}$.  We can
see this by comparing~(\ref{nov18.1}) and the first equation in~(\ref{dec1.2}).

\subsection{Two $\SLE_\kappa$ paths}

Consider the case of the $\SLE_\kappa$ measure corresponding
to two mutually avoiding $\SLE_\kappa$ paths; one from $0$
to $\infty$ and the other from $x_1$ to $x_2$.
This corresponds to choosing  
\begin{equation}  \label{dec24.1}
G(t,x_0,x_1,x_2)  
=  H^*(x_1-x_0,x_2-x_0)
= (x_2 - x_1)^{-2b} \, \phi_b\left(1-\frac{x_2 - x_1}{x_2 - x_0}\right)
\end{equation}
where $\phi_b$ is as in Proposition~\ref{dec1.prop1}. 
(Recall that $H^*(x,y)=H^*_{\Half,b,1}(x,y)$.) Then,
\[ 
\frac{\p_0 G(t,x_0,x_1,x_2)}
   {G (t,x_0,x_1,x_2)}   =
      - \frac{1}{x_2 - x_0}  \; \frac{q\, \phi_b'(1-q)}{\phi_b(1-q)}
  =  - \frac{q}{x_2 - x_0}\; [\log \phi_b]'(1-q), \]
where $\p_0=\p_{x_0}$ as before and $q = (x_2 - x_1)/(x_2 - x_0)$
denotes the appropriate cross-ratio.
 This is the marginal distribution of the $\SLE_\kappa$ path from
$0$ to $\infty$ given by the measure $\mu^\#_{b,2}(  (0,x_1),(x_2,\infty))$.
Two examples are as follows.
\begin{itemize}
\item
If $b=1$ (equivalently, $\kappa =2$), we have $\phi_1(u) = 2u-u^2$
and hence 
\[ \frac{\p_0 G(t,x_0,x_1,x_2)}
   {G (t,x_0,x_1,x_2)}   = -\frac{2}{x_2 - x_0}  \;
    \frac{q^2}{1-q^2} . \]
\item If $b = 1/4$ (equivalently, $\kappa = 4$), then $\phi_{1/4}(u) = \sqrt u$ and hence
\[   \frac{\p_0 G(t,x_0,x_1,x_2)}
   {G (t,x_0,x_1,x_2)}   = -\frac{1}{2(x_2 - x_0) }  \, \frac{q}{1-q}. \]
\end{itemize}

If we let $x_2 \rightarrow \infty$, we can consider two mutually avoiding
$\SLE_\kappa$ paths both going to infinity.  Here we use the 
function (see~(\ref{dec24.2}))
\[  G(t,x_0,x_1) = (x_1 - x_0)^a . \]
In this case,
\[ 
\frac{\p_0 G(t,x_0,x_1)}
   {G (t,x_0,x_1)}   = - \frac{a}{x_1 - x_0}. \]  From this we see that if
we take a pair mutually avoiding $\SLE_\kappa$ paths starting at
$x_0$, $x_1$, respectively, then the measure on the path started at $x_0$ is obtained 
by solving the usual chordal Loewner equation with driving function $U_t$
satisfying
\begin{equation*} 
      dU_t = \frac{a}{U_t - V_t} \, dt + dB_t,  
\end{equation*}
where $V_t = g_t(x_1)$ satisfies
\[           dV_t = \frac{a}{V_t - U_t}\,dt . \]
This is an example of a
 $\SLE(\kappa,\rho)$ process in~\cite{restriction}.
Note that $Z_t = U_t - V_t$ satisfies
\begin{equation}\label{apr28.eq1}
         dZ_t = \frac{2a}{Z_t} \,dt + dB_t. 
\end{equation}

\subsubsection{Growing both paths simultaneously}

It is possible to write down the equation for two mutually avoiding
$\SLE$ paths growing simultaneously.  We can write the corresponding
conformal transformation $g_t$ as the solution to the Loewner
equation
\[    \dot g_t(z) = \frac{a}{g_t(z) - U_t^1} + \frac{a}{g_t(z) - U_t^2}, \]
for appropriate driving functions $U_t^1$, $U_t^2$.  The parametrization
of $g_t$ is such that
\[  g_t(z) = z + \frac{2a}{z} + O(|z|^{-2}), \;\;\;\;\;  z \rightarrow \infty. \]
To find the equation of $U_t^j$ we start by taking them to satisfy
\[  dU_t^1 = \frac{a}{U_t^1 - U_t^2} \, dt + dW_t^1 \;\;\; \text{and} \;\;\;
  dU_t^2 = \frac{a}{U_t^2 - U_t^1} \, dt + dW_t^2, \]
where $W_t^1$, $W_t^2$ are independent standard Brownian motions.

Assume that $U_0^1 = x_0$ and $U_0^2 = x_1$.  If 
$G(t,x_0,x_1) = (x_1 - x_0)^a$,
 then
$G(t,U_t^1,U_t^2)$ satisfies
\begin{align*}
   \frac{dG(t,U_t^1,U_t^2)}{G(t,U^1_t,U_t^2)}  &=  
 \frac{a(a-1)}{(U_t^2 - U_t^1)^2} \, dt +  \frac{a}{U_t^1 - U_t^2}  \, dU_t^1
  + \frac{a} {U_t^2 - U_t^1}  \, dU_t^2  \\
 &=  \frac{a(3a -1)}{(U_t^2 - U_t^1)^2} \, dt
     +  \frac{a}{U_t^1 - U_t^2}  \, dW_t^1
  + \frac{a} {U_t^2 - U_t^1}  \, dW_t^2  
\end{align*}

Again,
\[  M_t = \exp\left\{-a(3a-1)\int_0^t (U_s^2 - U_s^1)^{-2}
\, ds\right\} \, G(t,U_t^1,U_t^2) \]
is a local martingale satisfying
\[  dM_t = M_t \,  \left[\frac{a}{U_t^1 - U_t^2}  \, dW_t^1
  + \frac{a} {U_t^2 - U_t^1}  \, dW_t^2\right]. \]
Compare this with~\cite[Example~1.18]{Lbook}. If we weight the paths by $M_t$, the effective equation on the $W_t^j$ is
\[  dW_t^j =    \frac{a} {U_t^j - U_t^{3-j}}   \, dt 
   + dB_t^j, \]
where $B_t^1$, $B_t^2$ denote independent standard Brownian motions, and the effective
equation on the $U_t^j$ becomes
\[       dU_t^j =    \frac{2a} {U_t^j - U_t^{3-j}}   \, dt 
   + dB_t^j. \]    
Note that if  
$Z_{t}= U_{t/2}^2 - U_{t/2}^1$,  $B_t = B_{t/2}^2 - B_{t/2}^1$, then
  $Z_t$, $B_t$ satisfy~(\ref{apr28.eq1}). 

\subsection{Multiple $\SLE_\kappa$}

Suppose $x_1  < \cdots <x_n < y_n < \cdots < y_1$.  Then the
measure $\mu^\#_{\Half,b,n}(\x,\y)$ can be described by using
the function
\[  G(t,x_1,\ldots,x_n,y_1,\ldots,y_n) =  H_{b,n}
  (\x,\y). \]
We can use this to describe the marginal distribution for
a particular path.  For example, suppose we are interested
in the first path.  Then it acts like $\SLE_\kappa$ weighted
locally by $G(t,W_t,X_t^2,\ldots,X_t^n,Y_t^1,\ldots,Y_t^n)$.


\section{The scaling limit of Fomin's identity}\label{fominscalesect}

The $\SLE$ measure is expected to be the limit of discrete
critical processes.  One example that is well understood at
this point is the case of the loop-erased random walk (LERW)
which is the prototypical model for the boundary scaling exponent
$b=1$.  We will describe
the discrete model in this section and
some different ways in which one can take the limit to get
the $\SLE$ measure with $b=1$ (corresponding to $\kappa=2$). 

In~\cite{fomin}, Fomin proved an identity showing that a 
certain probability for loop-erased random walk on $\Z^2$ can be 
given in terms of a determinant of hitting probabilities for simple 
random walk.  He also conjectured that this identity holds for continuous
 processes. In fact, he wrote,
 \begin{quote}
``\dots we do not need the notion of loop-erased Brownian motion. 
Instead, we discretize the model, compute the probability, and then pass to the limit."
\end{quote} 
It is shown in~\cite{KL} that the scaling limit of the determinant 
of hitting probabilities for simple random walk is, in fact, the determinant of 
the hitting densities of Brownian excursions, i.e., the so-called \emph{excursion 
Poisson kernel determinant}.

We conclude this section by deriving the scaling limit of Fomin's identity directly 
in the case $n=2$ using the ideas of the previous sections, and show how the probability 
that an $\SLE_2$ avoids a Brownian excursion may be computed in terms of the  
excursion Poisson kernel determinant. This is the natural continuous analogue of the
 discrete $n=2$ case, in which Fomin's identity gives the probability that 
a loop-erased walk avoids a simple random walk in terms of the determinant 
of the simple random walk hitting matrix.

\subsection{Loop-erased random walk}  \label{lerwsec}

We write the integer lattice $\Z^2$ in complex
form $\Z + i \Z$.  If $A \subsetneq \Z^2$, we let
$\p A =\{z \in \Z^2 : \dist(z,A) = 1\}$. We call
$\omega = [\omega_0,\ldots,\omega_n]$ a  {\em (nearest
neighbor)
random walk path} if $\omega_0,\ldots,\omega_n \in \Z^2$
with $|\omega_j - \omega_{j-1}| = 1$ for $j=1,\ldots, n$.
We write $|\omega| = n$ for the number of steps in the
path. If $z$, $w \in \p A$,
we say that $\omega$ is an excursion in $A$ from
$z$  to $w$ if $\omega_0 = z$, $\omega_n = w$, and
$\omega_1,\ldots,\omega_{n-1} \in A$. 
The {\em (simple random walk) excursion measure} $\excur(A;z,w)$
is the measure that assigns weight $4^{-|\omega|}$ to each
excursion $\omega$ from $z$ to $w$ in $A$.  We write
$h_A(z,w) = |\excur(A;z,w)|$.  The {\em loop-erased excursion}
measure $\excur_{LE}(A;z,w)$ is the measure obtained from
$\excur(A;z,w)$ by doing chronological loop-erasure on
each excursion. Note that $|\excur_{LE}(A;z,w)| =
 |\excur(A;z,w)| = h_A(z,w)$. We also write 
$\excur^\#(A;z,w) =  \excur(A;z,w)/h_A(z,w)$, and similarly for 
$\excur_{LE}^\#(A;z,w)$, to denote the measures normalized to be probability measures.

Let $\uplat = \{j + i k \in \Z^2: k > 0\}$ denote the
discrete upper half plane with $\p \uplat = \Z$. It is
not difficult to see (cf., equation~(28) in~\cite{KL}) that if $j_1 \neq j_2$, then
\[    h_\uplat(j_1,j_2) = \frac{1}{16} \, G_\uplat(j_1 + i, j_2 + i)\]
where $G_\uplat$ denotes the random walk Green's function in
$\uplat$.  Indeed we get the expression above by dividing
each excursion from $j_1$ to $j_2$ in $\uplat$ into its first
step $(j_1,j_1+i)$; its {\em middle}, i.e., any walk in $A$ starting
at $j_1 + i$ and ending at $j_2 + i$; and its final step $(j_2+i,j_2)$.
It is not hard to show (cf., equation~(15.3) in~\cite{Spitzer}) that
\[      G_\uplat(x,y) = a(x,\bar y) -  a(x,y) , \]
where $a(x,y) = a(y-x)$ denotes the potential kernel for random walk and $\bar y$
denotes the complex conjugate of $y$.  The potential kernel has
asymptotic behavior
\[              a(z) = \frac 2 \pi \, \log |z| + C + O(|z|^{-2}) , \]
for a known constant $C$.  Also, this relation can be differentiated, i.e., 
\[            a(z - 2i) - a(z) =  \frac 2 \pi
   \, [\log |z-2i| - \log|z|] + O(|z|^{-3}) . \]
Using this, we get the relation
\begin{equation}\label{LEexcmass}
   h_\uplat(j_1 , j_2) = \frac{1}{4 \pi \, (j_2-j_1)^2} +
   O((j_2 - j_1)^{-3}).
\end{equation}

Let us consider the scaling limit of this measure.   
  Suppose $N$ is a large integer.
   Let $\uplat_N = N^{-1} \uplat$ and let $\excur^{(N)}(j/N,k/N)$
be the measure that assigns measure $4\pi N^{2}4^{-|\omega|}$
to each excursion $\omega$ in $\uplat$ from $j$ to $k$. 
As $N \rightarrow \infty$, this measure approaches the {\em Brownian
excursion measure} $\excur_{B}$ on $\Half$ defined roughly
as follows.  The measure $\excur_{B}(x,y)$ has total mass
$H_1(x,y) = (y-x)^{-2}$. If $\excur_{B}^\#(x,y)$ denotes the
probability measure obtained by normalization, then this is the
measure corresponding to ``Brownian motion starting at $x$ conditioned
to leave $\Half$ at $y$'' which can be made precise using
$h$-processes. (See, for example,~\cite[Section~3.4]{Kozdron}.)

In the limit above, we fixed an $x$, $y$ and took the limit.  A similar
limit can be taken by considering the measure that gives
weight $4\pi4^{-|\omega|}  $
to every excursion (with any possible endpoints).  A similar argument
shows that the limiting measure obtained is (a constant multiple of)
Brownian excursion measure in $\Half$ and can be written as
\[          \int_{-\infty}^\infty \int_ {-\infty}
  ^\infty \frac{1}{(y-x)^2} \,
   \excur_{B}^\#(x,y) \; dx \, dy. \]
The advantage of this formulation is that it can be used to describe
the scaling limit of excursions in Jordan domains without any
smoothness assumptions on the boundary; see~\cite{Kozdron} and~\cite{KL}.

We can also take the scaling limit of the loop-erased measure
\[ \excur_{LE}(\uplat;j,k) = h_\uplat(j,k) \, \excur_{LE}^\#(\uplat;j,k). \]
The scaling limit of the normalized measure  $\excur_{LE}^\#(\uplat;j,k)$
was shown in~\cite{LSWlerw} to be $\SLE_2$ from $x$ to $y$.  (To be honest,
the paper~\cite{LSWlerw} considered a radial version, but the basic
argument  works  in this case.)

A {\em rooted (random walk) loop}
 in $A$ is a random walk path $\omega = [\omega_0,
\dots,\omega_{2n}]$ with $\omega_0 = \omega_{2n}$
and $\omega_1,\ldots,\omega_{2n} \in A$.   {\em An unrooted
(random walk)
loop} is an equivalence class of loops under the equivalence
\[  [\omega_0, \omega_1, \omega_2,
\dots,\omega_{2n}] \sim [\omega_1,\omega_2,\ldots,\omega_{2n}, \omega_1]. \]
For each unrooted loop $\omega$, let $R_\omega$ denote the number of distinct
rooted loops that produce this unrooted loop.  Note that $R_\omega$ is
always a divisor of $|\omega|$ and the proportion of loops of length
$2n$  with
$R_\omega < 2n$ decays exponentially in $n$. The random walk
loop measure in $A$ is the measure that assigns weight
$[R_\omega/|\omega|] \, 4^{-|\omega|}$ to each unrooted loop in $A$. 
If we scale the loops by $N^{-1}$ and take the limit, we get the
Brownian loop measure~\cite{rwloopsoup} as mentioned
in Section~\ref{Bloopsec}. 
 
The {\em random walk loop soup} in $\uplat$ is a Poisson point process
$\mathcal{C}_t$ from the loop measure.  At any fixed time $t$,
the set $\mathcal{C}_t$ is a multi-set of loops, where the number of
times the unrooted loop $\omega$ appears has a Poisson distribution
with parameter $(R_\omega/|\omega|) 4^{-|\omega|} t $. If $A
\subset \Z^2$, then ${\mathcal C}_t(A)$ denotes the loop soup in $A$,
i.e., ${\mathcal C}_t$ restricted to loops in $A$.
Given a self-avoiding path $\eta=[\eta_0,\eta_1,\ldots,\eta_k]$
and a realization of the loop soup ${ \mathcal C}_t$, there
is a systematic way of taking the loops in $\mathcal{C}_1(A)$ and adding
them to $\eta$ to produce a random walk path $\omega$ from
$\eta_0$ to $\eta_k$ whose chronological loop-erasure is $\eta$. 
It can be described as follows.
\begin{itemize}
\item Let $\omega^1,\ldots,\omega^r$ denote the loops in $\mathcal{C}_1(A)$
that intersect $\eta$.  It is not difficult to check that (with
probability one) there are
only a finite number of such loops (assuming $A$ is a proper
subset).  Order the loops in the order
that they appeared in the Poisson process.
\item  For each $\omega^j$ choose a rooted loop that represents it by
letting $\eta_l$ be the  point with smallest index $l$ such that $\eta_l$
appears in $\omega^j$ and choosing the root of $\omega$ so that the rooted path
starts at $\eta_l$.  If $\omega$ visits $\eta_l$ more than once,  choose
the root uniformly over all the possibilities. 
\item Attach a loop at $\eta_l$ by adding all the loops rooted at $\eta_l$
in the order they appear in the Poisson process.
\end{itemize}
Then straightforward
combinatorial arguments show that if $\eta$ is chosen according
to the probability measure $\excur_{LE}^\#(A;z,w)$ and $\mathcal{C}_1$ is an independent
realization of the loop soup in $A$ and the loops are added as above,
then the distribution of the resulting random walk is $\excur^\#(A;z,w)$.

Similarly, we can start with   $\excur^\#(A;z,w)$ and get
independent copies of the loop-erased measure  $\excur_{LE}^\#(A;z,w)$
and the loop soup in $A$ restricted to curves that hit the loop-erased
path.  We will explain this more now.  Suppose $\eta = [\eta_0,\eta_1,\ldots,
\eta_k]$ is a self-avoiding excursion in $A$ from $z$ to $w$.  The
$\excur_{LE}^\#(A;z,w)$ measure of $\eta$ is exactly the
$\excur^\#(A;z,w)$ measure of all loops $\omega$ whose loop-erasure
is $\eta$.   It is straightforward to show (see, for example, Equation~(12.2.2) of~\cite{LKesten}) that 
\[
\excur_{LE}^\#(A;z,w)(\eta) =4^{-|\eta|} \, \Theta_A(\eta) \;\;\; \text{ where } \;\;\;
         \Theta_A(\eta)^{-1}
   =    q_{A_1}(\eta_1) \cdots q_{A_{k-1}}(\eta_{k-1}).
\]
Here $A_j = A \setminus \{\eta_1,\ldots,\eta_{j-1}\}$ and $q_A(x)$
denotes the probability that a random walk starting at $x$ leaves $A$ before
returning to $x$.  We will show that 
\begin{equation}  \label{nov22.5}
 \Theta_A(\eta)  = {\Prob\{\mathcal{C}_1(A) \mbox{ contains no loop that
 intersects } \eta\}}.
\end{equation}

To verify~(\ref{nov22.5}), we
first note that any path $\omega$ with loop-erasure $\eta$ 
 can be split into a unique fashion as 
\[  \omega=[\eta_0,L_1,L_2,\ldots,L_{k-1},\eta_k] \]
where $L_j$ is a loop rooted at $\eta_j$ in $A \setminus \{\eta_1,
\ldots,\eta_{j-1}\}$.  Each of these loops can be split into
a number of ways as
\[           L_j = \omega^{j,1} \oplus \cdots \oplus \omega^{j,r_j} \]
where $\omega^{j,l}$ are loops at $L_j$.  We give $L_j$ weight
$4^{-|L_j|}$.  We let $\alpha( \omega^{j,l}) = 4^{-| \omega^{j,l}|}
  \theta( \omega^{j,l})^{-1}$ where $\theta([\omega_0,\ldots,\omega_{l}])$
is the number of $j \geq 1$ such that $\omega_j = \omega_0$.
We give the $l$-tuple $(\omega^{j,1},\ldots,\omega^{j,l})$
weight
\[  \alpha(\omega^{j,1}) \cdots \alpha(\omega^{j,l})
  = \frac 1{ 
4^{|\omega^{j,1}| + \cdots + |\omega^{j,r}| }\, [\theta(\omega^{j,1})
  \cdots \theta(\omega^{j,l})]} . \]
The following lemma yields~(\ref{nov22.5}).

\begin{lemma}
Suppose $A$ is a proper subset
of  $ \Z^2$ and $z \in A$.  Let $m=
m_{z,A}$ denote the measure
that assigns weight $4^{-|\omega|}\, \theta(\omega)^{-1}$ 
to each loop in $A$ rooted
at $z$.  Let $\mathcal{A}_t$ denote a Poisson point process from 
this measure and let $l_t$ denote the loop obtained
by concatenating all the loops in $\mathcal{A}_t$ in
chronological order. (If there are no loops, then $l_t$
is the trivial loop of length $0$.)  Then
\[  \Prob\{l_1 = \omega\} = 4^{-|\omega|} q_A, \]
where $q_A$ is the probability that simple random walk
starting at $z$ leaves $A$ before returning to $z$.
\end{lemma}

\begin{proof}
Let $\mathcal{L} = \mathcal{L}(A,z)$ denote the set of
all loops in $A$ rooted at $z$.   Note that
$   1 - q_A = 
           \sum_{\omega \in \mathcal{L}, \theta(\omega)
  = 1} 4^{-|\omega|} $
so that 
\begin{align}\label{dec17eq1} 
   \sum_{k=0}^\infty  
  \sum_{\omega \in \mathcal{L}, 
    \theta(\omega) = k } 4^{-|\omega|} 
   =   \sum_{k=0}^\infty 
   \left[\sum_{\omega \in \mathcal{L}, 
    \theta(\omega) = 1 } 4^{-|\omega|} \right]^k
 =  \sum_{k=0}^\infty  (1 - q_A)^k = 1/q_A. 
\end{align}
  For the trivial
loop $ [z]$,
$ \Prob\{l_1 = [z]\}$ is the probability that there
are no loops in the Poissonian realization, i.e.,
\[  \Prob\{l_1 = [z] \} = \exp\left\{
                 - \sum_{\omega \in \mathcal{L}} 4^{-|\omega|} 
   \theta(\omega)^{-1}
  \right\}.\]
It now follows from~(\ref{dec17eq1}) that
\begin{align*}
 \sum_{\omega \in \mathcal{L}} 4^{-|\omega|} 
   \theta(\omega)^{-1}
   =  \sum_{k=1}^\infty \frac 1k 
  \sum_{\omega \in \mathcal{L}, 
    \theta(\omega) = k } 4^{-|\omega|} 
& = \sum_{k=1}^\infty \frac 1k \, (1 - q_A)^k =
   -\log q_A 
\end{align*}
so therefore
\[      \Prob\{l_1 = [z] \} =  {q_A}.\]
More generally, if $\omega$ is a loop that returns to
the origin $k$ times, and $j_1,j_2, \ldots, j_r$
are positive integers summing to $k$, we could get
the loop $\omega$ by splitting the loop into $r$ loops
where the $m$th loop is from the $(j_1 + \cdots + j_{m-1})$th
to the $(j_1 + \cdots + j_m)$th return to $z$. By considering
all these possibilities, one can deduce that
\[  \Prob\{l_1 = \omega\} = q_A \, 4^{-|\omega|}
  \sum_{j_1 + \cdots +j_r
  = k} \frac{1}{r!} \, \frac{1}{j_1 \cdots j_r}. \]
We leave it to the reader to show that the summation equals
one and so $  \Prob\{l_1 = \omega\} = q_A \, 4^{-|\omega|}$.
(As a hint, write
\[    \frac{1}{1-t} = \exp\{-\log(1-t)\} \]
and expand both sides in power series in $t$.) 
\end{proof}

Computing the right side of~(\ref{nov22.5}) is very difficult;
the local terms corresponding to small loops dominate the probability.
However, we can immediately conclude a boundary perturbation
rule.  If $A \subset A'$ and $z,w \in \p A \cap \p A'$, then
\[  \frac{\Theta_{A'}(\eta)}
   {\Theta_A (\eta)} = \Prob\{\mathcal{C}_1(A') \mbox{ contains no loop 
intersecting both } A' \setminus A \mbox{ and } \eta \}. \]
This should be compared to Proposition~\ref{nov22.prop1}.

\subsection{Fomin's identity}

Suppose $A \subset \Z^2$ and
$x_1,\ldots,x_n,y_1,\ldots,y_n \in \p A$. Fomin~\cite{fomin} showed that the  
\[ \excur (A;x_1,y_1) \times \cdots \times \excur (A;x_n,y_n) \]
 measure  of the set of $n$-tuples of excursions
$(\omega^1,\ldots,\omega^n)$ satisfying the condition
\[            LE(\omega^{m-1}) \cap \omega^m =\emptyset,
\;\;\;\; m=2,\ldots,n, \]
where $LE(\omega^{m-1})$ denotes the loop erasure of $\omega^{m-1} $,
is given by 
 \[      \det[ h_{A}(x_l,y_m)]_{1 \leq l,m \leq n}. \]
In the case $A = \uplat$, topological considerations show that
this can be non-zero only if the points appear in order, i.e.,
after relabeling, they satisfy
\[   x_1 < \cdots <x_n < y_n < \cdots < y_1 . \]
As Fomin points out,  since the scaling limit of simple
random walk is Brownian motion, which is conformally invariant,
this formula in itself suggests the conformal invariance
of loop-erased random walk.  There are technical difficulties 
in taking the limit in domains
 with rough boundaries but these can
be handled; see~\cite{KL}. In the case $n=2$, the determinant equals
\[      h_{\uplat}(x_1,y_1) \, h_\uplat(x_2, y_2) \;
\left[ 1 
        - \frac{h_\uplat(x_1,y_2) \, h_\uplat(x_2,y_1)}
   { h_{\uplat}(x_1,y_1) \, h_\uplat(x_2, y_2)}\right] . \]
As the differences go to infinity, this is asymptotic to
\[         h_{\uplat}(x_1,y_1) \, h_\uplat(x_2, y_2) \;
    \left[1 -   
           \frac  
   {(y_1 - x_1)^2 \, (y_2 - x_2)^2}{(y_1-x_2)^2\, (y_2-x_1)^2}\right]. \]
If $x_1 = 0$, $x_2 = x$, $y_2 = y$,
 then as $y_1 \rightarrow \infty$ the last expression
is asymptotic to $h_{\uplat}(0,y_1)$ times
\[    h_{\uplat}(x,y) \, \frac{x}{y} \, \left[2 - \frac{x}{y}
  \right], \]
and we have recovered the expression for $H^*_{1,1}(x,y)$ that was derived
in Section~\ref{partsec}.

 \subsection{Review of excursion Poisson kernel determinant}
 
 We now briefly review the necessary notation from~\cite{KL} about the excursion Poisson kernel.  Suppose that $D \subset \C$ is a simply connected Jordan domain and that $\p D$ is locally analytic at $x$ and $y$. The \emph{excursion Poisson kernel} is defined as
 \[
 H_{\p D}(x,y) = \lim_{\eps \to 0} \frac{1}{\eps} \, H_D(x + \eps \n_x, y)
 \]
 where $H_D(z,y)$ for $z \in D$ is the usual Poisson kernel, and $\n_x$ is the unit normal at $x$ pointing into $D$. Explicit formul\ae\ are known when $D=\Disk$, the unit disk, or $D=\Half$, namely
 \[
 H_{\p \Disk}(x,y) = \frac{1}{\pi \, | y-x|^2} =\frac{1}{2\pi(1-\cos(\arg y - \arg x))}\;\;\; \text{ and } \;\;\; H_{\p \Half}(x,y) = \frac{1}{\pi(y-x)^2}.
\]
 Notice that $\pi H_{\p \Half}(x,y) = H_{\Half,1,1}(x,y)$. Since the excursion Poisson kernel is known~\cite[Proposition~2.11]{KL} to be conformally covariant, we obtain~(\ref{hinvariance}) with $b=1$, $n=1$.
Suppose now that $D' \subset \C$ is also a Jordan domain and that $f: D \to D'$ is a conformal transformation. Let $x_1, \ldots, x_n, y_1, \ldots, y_n$ be distinct boundary points at which $\p D$ is locally analytic, and assume that $\p D'$ is locally analytic at  $f(x_1), \ldots, f(x_n), f(y_1), \ldots, f(y_n)$.
It follows~\cite[Proposition~2.16]{KL} that if $\mathbf {H}_{\p D}(\x,\y) = [H_{\p D}(x_j, y_k)]_{1 \le j, k \le n}$ denotes the $n \times n$ \emph{hitting matrix}
\[
\mathbf {H}_{\p D}(\x,\y) = 
\begin{bmatrix}
H_{\p D}(x_1, y_1)   &\cdots    & H_{\p D}(x_1, y_n) \\
\vdots                         &\ddots    &\vdots \\
H_{\p D}(x_n, y_1)   &\cdots    &H_{\p D}(x_n, y_n) \\
\end{bmatrix}
 \]
 then
 \[ \det  \mathbf {H}_{\p D}(\x,\y) = \left( \prod_{\ell=1}^n |f'(x_\ell)| \; |f'(y_\ell)| \right) \det  [H_{\p D'}(f(x_j), f(y_k))]_{1 \le j, k \le n}. \]
 
\subsection{Direct derivation of the scaling limit}\label{directderivsect}

Suppose that $\gamma$ is an $\SLE_2$ in $\Half$ from $0$ to $\infty$ as in Section~\ref{SLEdefnsect}, and for every $0 < t < \infty$ let $\Half_t$ denote the slit-plane
$\Half_t = \Half \setminus \gamma[0,t]$.  Fix two real numbers $0<x<y<\infty$, and let $\beta:[0,1] \to \overline{\Half}$ be a Brownian excursion from $x$ to $y$ in $\Half$. It then follows that
\begin{equation}\label{probtocompute}
\Prob\{ \,\gamma[0,t] \cap \beta[0,1] = \emptyset \; | \; \gamma[0,t] \,\} =\frac{H_{\p \Half_t}(x,y)}{H_{\p \Half}(x,y)} = \frac{H_{\Half_t, 1,1}(x,y)}{H_{\Half, 1,1}(x,y)}
\end{equation}
so that taking expectations and limits of~(\ref{probtocompute}) gives 
 \begin{equation}\label{fominprob}
 \Prob\{ \,\gamma[0,\infty) \cap \beta[0,1] = \emptyset\, \} = \tilde H^*_{\Half,1,1}(x,y) = \frac{x}{y} \left(2 - \frac{x}{y} \right).
\end{equation}
 The first equality is similar to~(\ref{dec16eq1}), and the second equality follows from Proposition~\ref{dec1.prop1}.
 
 \begin{theorem}
If $x$, $y \in \R$ with $0<x<y<\infty$, and that $\gamma$, $\beta$ are as above, then \begin{equation}\label{scalefomin}
 \Prob\{ \,\gamma[0,\infty) \cap \beta[0,1] = \emptyset\, \}
=
\frac{\det  \mathbf {H}_{\p \Disk}(f(\x),f(\y)) }{ 
H_{\p \Disk}(f(0), f(\infty))H_{\p \Disk}(f(x), f(y))
}
\end{equation}
where $f: \Half \to \Disk$ is a conformal transformation.
\end{theorem}
 
 \begin{remark}
 By working in $\Half$ and $\Disk$, it is possible to perform explicit calculations. Since the quantity on the right side of~(\ref{scalefomin}) is known to be a conformal invariant, we can show, with a combination of conformal transformations, that the probability that an $\SLE_2$ avoids a Brownian excursion in $D$ is given by the appropriate determinant of the matrix of excursion Poisson kernels.
 \end{remark}

\begin{proof}
By the scaling rule~(\ref{scalerule}), it suffices without loss of generality to assume that $0<x<1$ and $y=1$. Furthermore,
we may assume that the conformal transformation $f: \Half \to \Disk$ is given by
\[f(z) = \frac{iz+1}{z+i},\]
so that $f(0)=-i$, $f(1)=1$, $f(\infty)=i$, and
\[
f(x) = \left(\frac{2x}{x^2+1} \right) + i \left( \frac{x^2-1}{x^2+1} \right) = \exp \left\{ -i \arctan \left( \frac{1-x^2}{2x} \right) \right\}.
\]
Notice that $(2x)^2+(x^2-1)^2=(x^2+1)^2$ so that $|f(x)|=1$ as expected. Writing $f(x) = e^{i \theta}$, we find that
\begin{align*}
\frac{\det  \mathbf {H}_{\p \Disk}(f(\x),f(\y)) }{ 
H_{\p \Disk}(f(0), f(\infty))H_{\p \Disk}(f(x), f(y)) }
&= \frac{\det  \mathbf {H}_{\p \Disk}(-i,i,f(x),1) }{ 
H_{\p \Disk}(-i,i)H_{\p \Disk}(f(x), 1) }\\
&= \frac{
H_{\p \Disk}(-i,i)H_{\p \Disk}(e^{i\theta},1)-H_{\p \Disk}(-i,1)H_{\p \Disk}(e^{i\theta},i)}
{H_{\p \Disk}(-i,i)H_{\p \Disk}(e^{i\theta},1) }\\
&= \frac{ \frac{1}{2\pi (1-\cos \pi)} \frac{1}{2\pi (1-\cos \theta)}
-\frac{1}{2\pi (1-\cos (\frac{\pi}{2}) )} \frac{1}{2\pi (1-\cos(\frac{\pi}{2}+\theta))}}
{\frac{1}{2\pi (1-\cos \pi)} \frac{1}{2\pi (1-\cos \theta)}}\\
&= \frac{2\cos \theta + \sin \theta - 1}{1 + \sin \theta}.
\end{align*}
Since 
$\theta = - \arctan \left(\frac{1-x^2}{2x} \right)$
we see that
$\cos \theta  = \frac{2x}{x^2+1}$  and $\sin \theta = \frac{1-x^2}{x^2+1}$
which upon substitution gives
\[ 
\frac{2\cos \theta + \sin \theta - 1}{1 + \sin \theta}
= \frac{ \frac{4x}{x^2+1} + \frac{1-x^2}{x^2+1} - 1}{1 +  \frac{1-x^2}{x^2+1}}
= \frac{4x-2x^2}{2} = x(2-x).\]
Comparison with~(\ref{fominprob}) now yields the result.
\end{proof}

\section{The $\lambda$-SAW}

Here we will define a one-parameter family of measures on self-avoiding
walks on the
integer lattice
$\Z^2$.  We conjecture that the scaling limit of this measure will
give the measure on $\SLE$ as described before.  

A simple random walk path $\omega$ that does not visit any point more
than once is called a self-avoiding walk.  If $\omega$ is an excursion
in $A$ connecting boundary points $z$, $w$, we will call $\omega$
a self-avoiding excursion (SAE).  Let $r$, $\lambda \in \R$ be parameters.
If $A \subsetneq \Z^2$ and $z$, $w \in \p A$,
then the measure $\msae_{r,\lambda,z,w,A}$ is defined to be the
measure that assigns to each SAE $\omega$
 in $A$ from $z$ to $w$ measure
\[     \exp\{-r|\omega| + \lambda m^*(A;\omega)\} , \]
where $m^*(A;\omega)$ denotes the random walk loop measure of the set
of loops in $A$ that intersect $\omega$.  In analogy
with the discussion on $\SLE$, we conjecture that for each
$\lambda \geq -1/2$, there is a critical $r_\lambda$ such that  
  the total mass
$|\msae_{r_\lambda,\lambda,0,N,\uplat}|$ neither decays nor
grows exponentially as $N \rightarrow \infty$.  
\begin{itemize}
\item  If $\lambda = 1$ and $e^{-r} = 1/4$,
then   $\msae_{r,1,0,N,\uplat}$ is exactly the same as
the loop-erased excursion measure 
$\excur_{LE}(\uplat;0,N)$
which from~(\ref{LEexcmass}) has total mass asymptotic
to $ (4\pi)^{-1}  N^{-2}$.  Therefore, $e^{-r_1} = 1/4$.
\item  If $\lambda = 0$, then this measure is the same as the
usual self-avoiding walk measure that gives each walk of $n$ steps
the same measure.  In this case $e^{r_0}$ is the {\em connective
constant} for self-avoiding walk which can be defined by saying that the
number of SAWs of length $n$ starting at the origin in $\Z^2$ grows
like $e^{r_0n}$. 
\end{itemize}

We now make this conjecture stronger and say that there exist $ C >0$, 
$b \in \R$ (depending on $\lambda$) such that
$|\msae_{r_\lambda,\lambda,0,N,\uplat}| \sim C \, N^{-2b}$.
Consider the measure $N^{2b-2} \, C^{-1} \, \msae_{r_\lambda,\lambda,
\cdot,\cdot,\uplat}$ where we use $\cdot,\cdot$ to mean that we take
the union over all $x$, $y \in \R$.  We also scale the paths by $N^{-1}$.
Then we conjecture that we have a measure on excursions
in $\Half$.  If $x_1 < x_2 < x_3 < x_4$, then the total mass of
the set of excursions connecting $[x_1,x_2]$ to $[x_3,x_4]$ should be
\[             \int_{x_1}^{x_2} \int_{x_3}^{x_4}\, \frac{1}{(t-s)^{2b}}
          dt \, ds . \]
We conjecture that we can write this measure as
\[             \int_{x_1}^{x_2} \int_{x_3}^{x_4}\, \frac{1}{(t-s)^{2b}}
\, \mu^\#_\Half(s,t)\; 
          dt \, ds , \]
where $\mu^\#_\Half(s,t)$ is a conformally invariant probability measure
on paths.

We have not discussed the parametrization of the limiting measure.  An
even stronger conjecture would say that there is an exponent\footnote{It is traditional to use $\nu$ for this critical exponent. Although we have also used $\nu$ with parameters to denote certain measures on paths, this should not cause confusion.} $\nu$ defined
roughly by saying that the typical number of steps in a walk from $0$
to $N$ in the measure above contains $N^{1/\nu}$ steps. Then we could
let
\[   \omega^{(N)}(t) = N^{-t} \omega(t \, N^{1/\nu}). \]
However, it is not necessary to include the parametrization to have
a nontrivial result.

The measure on $k$-tuples of self-avoiding walks is defined similarly,
by giving measure 
\[        \exp\{-r(|\omega_1| + \cdots + |\omega_k|)  
 + \lambda m^*(A;\omega_1 \cup \cdots \cup \omega_k)\} . \]

We conjecture the following.
\begin{itemize}
\item  If $\lambda \geq -1/2$, then the $\lambda$-SAW gives a limiting
measure on paths that is the $\SLE_\kappa$ measure.  In this
case
$\kappa,\lambda,b$ are related as before.
\item  For $\lambda < -1/2$, this does not give a limiting measure
on simple paths.
\end{itemize}
For $\lambda = 1$, this has been proved using the loop-erased random
walk.  For $\lambda = 0$, this is equivalent to very difficult conjectures
about self-avoiding walks; see~\cite{LSWsaw}.


\section*{Acknowledgements}
The authors would like to express their gratitude to the Fields Institute for their gracious hospitality during the \emph{Percolation, $\SLE$, and Related Topics Workshop} of September~2005.


\end{document}